\documentclass[reqno,10pt]{amsart}
\usepackage{amsmath,amssymb}

\newtheorem{thm}{Theorem}[section]
\newtheorem{lem}[thm]{Lemma}
\newtheorem{cor}[thm]{Corollary}
\newtheorem{defi}[thm]{Definition}
\theoremstyle{definition}
\newtheorem{rek}[thm]{Remark}

\newcommand\ben{\begin{enumerate}}
\newcommand\een{\end{enumerate}}

\newcommand{\glfour}{{\rm GL}(4)}
\newcommand{\glsix}{{\rm GL}(6)}

\newcommand{\soe}{{\rm SO(even)}}
\newcommand{\soo}{{\rm SO(odd)}}
\newcommand{\sym}{{\rm sym}}
\newcommand{\SO}{{\rm SO}} 

\newcommand{\R}{\mathbb{R}}

\newcommand{\ga}{\alpha}    
\newcommand{\gb}{\beta}      
\newcommand{\gd}{\delta}    
\newcommand{\gep}{\epsilon}  
\newcommand{\g}{\gamma}      
\newcommand{\gl}{\lambda}    
\newcommand{\gt}{\theta}    

\newcommand{\twocase}[5]{#1 \begin{cases} #2 & \text{\rm #3}\\ #4
&\text{\rm #5} \end{cases}  }

\newcommand\be{\begin{equation}}
\newcommand\ee{\end{equation}}
\newcommand\bea{\begin{eqnarray}}
\newcommand\eea{\end{eqnarray}}

\newcommand{\smallsubstack}[2]{\genfrac{}{}{0pt}{}{#1}{#2}}

\newcommand{\foh}{\frac{1}{2}}  

\newcommand{\hkpn}{H_k^+(N)}
\newcommand{\hkmn}{H_k^-(N)}
\newcommand{\hkpmn}{H_k^\pm(N)}
\newcommand{\hkn}{H_k^\ast(N)}
\newcommand{\hksn}{H_k^\sigma(N)}

\newcommand{\jk}[2]{J_{k-1}\left( 4\pi \frac{ \sqrt{ #1 } }{ #2 }\right) }
\newcommand{\phir}[1]{\widehat{\phi}\left( \frac{ \log p_{#1} }{\log R}\right) }

\newcommand{\pfrac}[1]{\frac{2\log p_{#1}}{\sqrt{p_{#1}} \log R}}

\newcommand{\ils}{\cite{ILS}}

\newcommand{\db}{\overline{d}}

\newcommand{\hphi}{\widehat{\phi}}  
\newcommand{\Zf}{Z_\phi} 

\renewcommand{\mod}{\;\operatorname{mod}}
\newcommand{\smod}[1]{(\operatorname{mod} #1)}
\newcommand{\supp}{\operatorname{supp}}
\newcommand{\notdiv}{\nmid}
\newcommand{\intinf}{\int_{-\infty}^\infty}
\newcommand{\E}{{\mathbb E}} 
\newcommand{\I}{1\!\!1} 

\renewcommand{\i}{{\mathrm{i}}} 
\renewcommand{\d}{{\mathrm{d}}} 

\renewcommand{\Re}{{\mathfrak{Re}}}
\renewcommand{\Im}{{\mathfrak{Im}}}

\newcommand{\<}{\left\langle}
\renewcommand{\>}{\right\rangle}

\numberwithin{equation}{section}

\begin{document}

\title{Low lying zeros of $L$-functions with orthogonal symmetry}

\author{C.P. Hughes}
\address{Department of Mathematics, University of Michigan, Ann
Arbor, MI 48109}
 \email{hughes@aimath.org}

\author{Steven J. Miller}
\address{Department of Mathematics, Brown University, 151 Thayer
 Street,
Providence, RI 02912}
 \email{sjmiller@math.brown.edu}



\date{\today}

\begin{abstract}
We investigate the moments of a smooth counting function of the
zeros near the central point of $L$-functions of weight $k$ cuspidal
newforms of prime level $N$. We split by the sign of the functional
equations and show that for test functions whose Fourier transform
is supported in $(-\frac1{n},\frac1{n})$, as $N\to\infty$ the first
$n$ centered moments are Gaussian. By extending the support to
$(-\frac{1}{n-1},\frac1{n-1})$, we see non-Gaussian behavior; in
particular the odd centered moments are non-zero for such test
functions. If we do not split by sign, we obtain Gaussian behavior
for support in $(-\frac{2}{n}, \frac{2}{n})$ if $2k \ge n$.  The
$n$\textsuperscript{th} centered moments agree with Random Matrix
Theory in this extended range, providing additional support for the
Katz-Sarnak conjectures. The proof requires calculating
multidimensional integrals of the non-diagonal terms in the
Bessel-Kloosterman expansion of the Petersson formula. We convert
these multidimensional integrals to one-dimensional integrals
already considered in the work of Iwaniec-Luo-Sarnak, and derive a
new and more tractable expression for the $n$\textsuperscript{th}
centered moments for such test functions. This new formula
facilitates comparisons between number theory and random matrix
theory for test functions supported in $(-\frac1{n-1},\frac1{n-1})$
by simplifying the combinatorial arguments. As an application we
obtain bounds for the percentage of such cusp forms with a given
order of vanishing at the central point.
\end{abstract}

\maketitle

\section{Introduction}\label{sect:introduction}

Let $H^\star_k(N)$ be the set of all holomorphic cusp forms of
weight $k$ which are newforms of level $N$. Every $f\in
H^\star_k(N)$ has a Fourier expansion
\begin{equation}
f(z)\ =\ \sum_{n=1}^\infty a_f(n) e(nz).
\end{equation}
Set $\lambda_f(n) =  a_f(n) n^{-(k-1)/2}$. The $L$-function
associated to $f$ is
\begin{equation}
L(s,f)\ =\ \sum_{n=1}^\infty \lambda_f(n) n^{-s}.
\end{equation}
The completed $L$-function is
\begin{equation}\label{eq:completed_L_func}
\Lambda(s,f) \ =\ \left(\frac{\sqrt{N}}{2\pi}\right)^s
\Gamma\left(s+\frac{k-1}{2}\right) L(s,f),
\end{equation}
and it satisfies the functional equation $\Lambda(s,f) =
\epsilon_f \Lambda(1-s,f)$ with $\epsilon_f = \pm 1$. Therefore
$H^\star_k(N)$ splits into two disjoint subsets, $H^+_k(N) = \{
f\in H^\star_k(N): \epsilon_f = +1\}$ and $H^-_k(N) = \{ f\in
H^\star_k(N): \epsilon_f = -1\}$. Each $L$-function has a set of
non-trivial zeros $\rho_f = \tfrac12 + \i\g_f$. The Generalized
Riemann Hypothesis is the statement that all $\g_f \in \R$ for all
$f$.

Assuming GRH, the zeros of any such $L$-function lie on the critical
line, and therefore it is possible to investigate statistics of the
normalized zeros (that is, the zeros which have been stretched out
to be one apart on average). The general philosophy, born out in
many examples (see for example \cite{CFKRS,KeSn}), is that
statistical behavior of families of $L$-functions can be modeled by
ensembles of random matrices. The spacing statistics of high zeros
of automorphic cuspidal $L$-functions (see \cite{Mon,Hej,RS}), for
certain test functions, agree with the corresponding statistics of
eigenvalues of unitary matrices chosen with Haar measure (or
equivalently, complex Hermitian matrices whose independent entries
are chosen according to Gaussian distributions). Initially this led
to the belief that this was the only matrix ensemble relevant to
number theory; however Katz and Sarnak (\cite{KS1,KS2}) prove that
these statistics are the same for all classical compact groups.
These statistics, the $n$-level correlations, are insensitive to
finitely many zeros; thus, differences in behavior at the central
point $s = \foh$ are missed by such investigations, and a new
statistic, sensitive to behavior near the central point, is needed
to distinguish families of $L$-functions. In many cases
(\cite{ILS,Ru,Ro,HR2,FI,Mil,Yo,DM,Gu,Gao}) the behavior of the low
lying zeros (zeros near the central point) of families of
$L$-functions are shown to behave similarly to eigenvalues near $1$
of classical compact groups (unitary, symplectic and orthogonal).
The different groups exhibit different behavior near $1$.

Let $\phi$ be an even Schwartz function such that its Fourier
transform has compact support. We are interested in moments of the
smooth counting function (also called the one-level density or
linear statistic)
\begin{equation}\label{eq:dfphidef}
D(f;\phi)\ =\ \sum_{\g_f} \phi\left(\frac{\log R}{2\pi}
\g_f\right)
\end{equation}
when averaged over either $H_k^\ast(N)$ (the unsplit case), or
$H_k^+(N)$ or $H_k^-(N)$ (the split cases) as $N\to\infty$ through
the primes, with $k$ held fixed. Here $\gamma_f$ runs through the
non-trivial zeros of $L(s,f)$, and $R$ is its analytic conductor
($R = k^2N$ for these families). We rescale the zeros by $\log R$
as this is the order of the number of zeros with imaginary part
less than a large absolute constant. Because of the rapid decay of
$\phi$, most of the contribution in \eqref{eq:dfphidef} is from
zeros near the central point. We
use the uniform average over $f\in H_k^\sigma(N)$ (for
$\sigma$ one of $\ast$, $+$ or $-$), in the sense that if $Q$ is a function defined
on $f\in H_k^\sigma(N)$, then the average of $Q$ over
$H_k^\sigma(N)$ is
\begin{equation}\label{eq:unifaverdef}
\< Q(f) \>_\sigma\ :=\ \frac{1}{|H_k^{\sigma}(N)|} \sum_{f\in
H_k^\sigma(N)} Q(f).
\end{equation} We discuss in detail in Remarks
\ref{rek:whynoharm} and \ref{rek:HoffLocknoharm} why we chose to
uniformly weigh each $f \in H_k^\sigma(N)$ and not use harmonic
averaging as in \cite{Ro}, though both approaches yield the same
support.

Our first theorem evaluates the centered moments of $D(f,\phi)$
over $f \in H_k^\ast(N)$.

\begin{thm}\label{thm:nosplitmoments}
Assume GRH for $L(s,f)$. For $n\geq 1$ an integer, if $\supp(\hphi) \subset
(-\frac1{n}\frac{2k-1}{k},\frac1{n}\frac{2k-1}{k})$ then the $n${\rm \textsuperscript{th}}
centered moment of $D(f;\phi)$, when averaged over the elements of
$H_k^\ast(N)$, converges as $N\to\infty$ through
prime values to the $n${\rm \textsuperscript{th}} centered moment of
the Gaussian distribution with mean
\begin{equation}
\hphi(0)+\frac12\int_{-1}^1 \hphi(y)\;\d y
\end{equation}
and variance
\begin{equation} \sigma_\phi^2 \ = \
2\int_{-1/2}^{1/2} |y| \hphi(y)^2\;\d y.
\end{equation}
\end{thm}

\begin{rek}\label{rek:whyGRHforLsf}
We assume GRH for $L(s,f)$ to simplify the arguments below; however,
we may remove this assumption by arguing as on page 88 of \cite{ILS}
(specifically, either use the Petersson formula to handle the $p^2$
terms in the explicit formula, or crude estimates for $L(s, \sym^2 f
\otimes \sym^2 f)$). \end{rek}

Thus, with a little more work, Theorem \ref{thm:nosplitmoments} can
be made unconditional. By assuming GRH for Dirichlet $L$-functions,
in Theorem \ref{thm:nosplitbetterk} we increase the support to
$(-\frac2n, \frac2n)$, provided $2k \ge n$. The relation between $n$
and $k$ arises from some technicalities in controlling error terms;
these obstructions are usually not apparent in studying just the
$n=1$ case.

If we split by sign, then the same argument still gives Gaussian
moments, but with a greater restriction on the support of the test
function $\hphi$. Later we increase the support by invoking GRH for
Dirichlet $L$-functions.

\begin{thm}\label{thm:mock-Gaussian for D}
Under GRH for $L(s,f)$, if $\supp(\hphi) \subset
(-\frac1{n},\frac1{n})$ then the $n${\rm \textsuperscript{th}}
centered moment of $D(f;\phi)$, when averaged over the elements of
either $H^+_k(N)$ or $H_k^-(N)$, converges as $N\to\infty$ through
prime values to the $n${\rm \textsuperscript{th}} centered moment of
the Gaussian distribution with mean
\begin{equation}
\hphi(0)+\frac12\int_{-1}^1 \hphi(y)\;\d y
\end{equation}
and variance
\begin{equation} \sigma_\phi^2 \ = \
2\int_{-1/2}^{1/2} |y| \hphi(y)^2\;\d y.
\end{equation}
\end{thm}

Hughes and Rudnick \cite{HR1} prove a similar result within random
matrix theory. For a Schwartz function $\phi$ on the real line,
define
\begin{equation}
F_M(\theta)\ := \ \sum_{j=-\infty}^\infty \phi\left(\frac
M{2\pi}(\theta+2\pi j)\right),
\end{equation}
which is $2\pi$-periodic and localized on a scale of $\frac1{M}$.
For $U$ an $M\times M$ unitary matrix with eigenvalues
$e^{\i\theta_n}$, set
\begin{equation}
\Zf(U)\ :=\ \sum_{n=1}^M F_M(\theta_n).
\end{equation}
Note that going from $e^{\i\theta_n}$ to $\theta_n$ is well
defined, since $F_M(\theta)$ is $2\pi$-periodic. We often consider
$U$ to be a special orthogonal matrix when the eigenvalues occur
in complex-conjugate pairs, and thus will be doubly counted.
$\Zf(U)$ is the random matrix equivalent of $D(f;\phi)$. In
\cite{HR1} it was proved that
\begin{thm}\label{thm:mock_Gaussian}
If $\supp(\hphi) \subseteq [-\frac1{n},\frac1{n}]$ then the first
$n$ moments of $\Zf(U)$ averaged with Haar measure over $\SO(M)$
(with $M$ either even or odd) converge to the moments of the
Gaussian distribution with mean
\begin{equation}\label{eq:RMT_mean}
\hphi(0)+\frac12\int_{-1}^1 \hphi(y)\;\d y
\end{equation}
and variance
\begin{equation}\label{eq:RMT_var}
\sigma_\phi^2 \ = \
2\int_{-1/2}^{1/2} |y| \hphi(y)^2\;\d y.
\end{equation} In particular, the odd moments vanish, and for $2m \le n$ the
$2m^{\rm th}$ moment is $(2m-1)!!\ \sigma^{2m}_\phi$.

If $\supp \hphi \subseteq \left[-\frac{2}{n}, \frac{2}{n}\right]$,
then the $n^{\rm th}$ moment of $\Zf(U)$ when averaged over the
mean\footnote{By the mean of $\soe$ and $\soo$ we mean the ensemble
where half the matrices are $\soe$ and the other half $\soo$.} of
$\soe$ and $\soo$ converges to the $n^{\rm th}$ moment of a Gaussian
random variable with mean and variance given by \eqref{eq:RMT_mean}
and \eqref{eq:RMT_var} respectively.
\end{thm}

Thus Theorems \ref{thm:nosplitmoments} (and
\ref{thm:nosplitbetterk}), \ref{thm:mock-Gaussian for D} and
\ref{thm:mock_Gaussian} provide evidence for the connection between
number theory and random matrix theory, specifically that the
behavior of zeros near the central point is well modeled by that of
eigenvalues near $1$ of a classical compact group. It was remarked
in \cite{HR1} that the $n$\textsuperscript{th} moment of $\Zf$ when
averaged over either $\soe$ or $\soo$ ceases to be Gaussian once the
support of $\hphi$ is greater than $[-\frac{1}{n},\frac{1}{n}]$ (we
make this remark precise in Theorem~\ref{thm:mmts_RMT}). Similarly
we prove Theorem \ref{thm:mock-Gaussian for D} is essentially sharp
by showing the odd centered moments of $D(f;\phi)$ are non-zero if
the support of $\hphi$ is strictly greater than
$[-\frac{1}{n},\frac{1}{n}]$, and thus the distribution of
$D(f;\phi)$ is non-Gaussian in that case. Furthermore, the
$n$\textsuperscript{th} centered moments when averaged over
$H_k^+(N)$ or $H_k^-(N)$ are different as soon as the support of
$\hphi$ exceeds $[-\frac1n, \frac1n]$. This was anticipated in the
work of Iwaniec, Luo and Sarnak \cite{ILS}, who proved the following
theorem (Theorem 1.1 of \cite{ILS}):

\begin{thm}\label{thm:mean of D}
If $\supp(\hphi) \subset (-2,2)$, then the first moment agrees with
random matrix theory; explicitly, under GRH for $L(s,f)$ and for
Dirichlet $L$-functions,
\begin{align}\label{eq:mean of D}
\lim_{N\to\infty} \< D(f;\phi) \>_+ &\ = \ \hphi(0) +
\frac12\int_{-1}^1 \hphi(y) \;\d
y \\
\lim_{N\to\infty} \< D(f;\phi) \>_- &\ = \ \hphi(0) +
\frac12\int_{-1}^1 \hphi(y) \;\d y + \int_{|y| \geq 1}\hphi(y)
\;\d y.
\end{align}
(Here $N$ tends to infinity through the square-free integers).
\end{thm}

While the assumption of GRH for Dirichlet $L$-functions is essential
above (it is needed to analyze the Kloosterman sums), on page 88 of
\cite{ILS} they give two arguments which remove the assumption of
GRH for $L(s,f)$.

Note that if $\supp(\hphi) \subset (-1,1)$ then $\lim \< D(f;\phi)
\>_+ = \lim \< D(f;\phi) \>_-$, but they are different if the
support of $\hphi$ is outside this interval. Thus in order to test
the expected belief that averages over $H_k^+(N)$ correspond to
averages over $\soe$ and averages over $H_k^-(N)$ correspond to
averages over $\soo$, it is essential that the calculations in
\cite{ILS} have support greater than $1$, as for smaller support the
$1$-level densities of the orthogonal groups are indistinguishable.
In this paper we further test this correspondence. Our main result
is
\begin{thm}\label{thm:extended moments}
Let $n \geq 2$, $\supp(\hphi) \subset (-\frac1{n-1},
\frac1{n-1})$, $D(f;\phi)$ be as in \eqref{eq:dfphidef}, and
define
\begin{align}\label{eq:defRnn-1}
R_n(\phi) \ &= \ (-1)^{n-1} 2^{n-1}\left[ \int_{-\infty}^\infty
\phi(x)^{n}
\frac{\sin 2\pi x}{2\pi x} \;\d x - \foh \phi(0)^{n} \right]\\
\sigma^2_\phi\label{eq:defRnn-1b}  \ &= \ 2\int_{-1}^{1} |y|
\hphi(y)^2\;\d y.
\end{align}
Assume GRH for $L(s,f)$ and for all Dirichlet $L$-functions. As
$N\to\infty$ through the primes,
\begin{equation}\label{eq:thmextmomcompact}
\lim_{\substack{N\to\infty \\ N \text{\rm prime}}} \<
\left(D(f;\phi) - \< D(f;\phi) \>_\pm \right)^{n}\>_\pm  =
\begin{cases}
(2m-1)!!\ \sigma^{2m}_\phi \pm R_{2m}(\phi) &
\text{\rm if $n=2m$ is even}\\
\pm R_{2m+1}(\phi) & \text{\rm if $n=2m+1$ is odd.}
\end{cases}
\end{equation}
\end{thm}

Note if $\supp(\hphi) \subset (-\frac1n,\frac1n)$ then $R_n(\phi) =
0$ and we recover the Gaussian behavior of
Theorem~\ref{thm:mock-Gaussian for D}. Also $R_n(\phi)$ is not
identically zero for test functions $\phi$ such that $\supp(\hphi)
\not\subset [-\frac1n,\frac1n]$. The assumption of GRH for $L(s,f)$
is for ease of exposition, and can be removed; see Remark
\ref{rek:whyGRHforLsf}.

Finally we show that the random matrix moments of $\Zf$ correctly
model the moments of $D(f;\phi)$ (at least for the support of
$\hphi$ restricted as in Theorem~\ref{thm:extended moments}), in the
sense that the $n$\textsuperscript{th} centered moment of
$D(f;\phi)$ averaged over $H_k^+(N)$ equals the
$n$\textsuperscript{th} centered moment of $\Zf$ averaged over
$\soe$, and $H_k^-(N)$ similarly corresponds to $\soo$.

\begin{thm}\label{thm:mmts_RMT} The means of
$\Zf(U)$ when averaged with respect to Haar measure over $\soe$ or
$\soo$ are
\begin{align}
\mu_+ &:= \lim_{\substack{M\to\infty \\ M {\rm even}}} \E_{\SO(M)}
\left[ \Zf(u) \right]\ =\
\hphi(0)+\frac12\int_{-\infty}^\infty \hphi(y)\;\d y \\
\mu_- &:= \lim_{\substack{M\to\infty \\ M {\rm odd}}} \E_{\SO(M)}
\left[ \Zf(u) \right]\ = \ \hphi(0)+\frac12\int_{-\infty}^\infty
\hphi(y)\;\d y + \int_{|y|\geq 1} \hphi(y) \;\d y.
\end{align}
Let $R_n(\phi)$ and $\sigma^2_\phi$ be as in \eqref{eq:defRnn-1} and
\eqref{eq:defRnn-1b}. For $n \geq 2$ if $\supp(\hphi) \subseteq
[-\frac{1}{n-1} , \frac{1}{n-1}]$ then the $n${\rm
\textsuperscript{th}} centered moment of $\Zf(U)$ converges to
\begin{equation}
\lim_{\substack{M\to\infty \\ M {\rm even}}} \E_{\SO(M)} \left[
\left(\Zf(U) - \mu_+ \right)^n \right] =
\begin{cases}
(2m-1)!!\ \sigma^{2m}_\phi + R_{2m}(\phi) &
\text{\rm if $n=2m$ is even}\\
R_{2m+1}(\phi) & \text{\rm if $n=2m+1$ is odd}
\end{cases}
\end{equation}
and
\begin{equation}
\lim_{\substack{M\to\infty \\ M {\rm odd}}} \E_{\SO(M)} \left[
\left(\Zf(U) - \mu_- \right)^n \right] =
\begin{cases}
(2m-1)!!\ \sigma^{2m}_\phi - R_{2m}(\phi) &
\text{\rm if $n=2m$ is even}\\
-R_{2m+1}(\phi) & \text{\rm if $n=2m+1$ is odd.}
\end{cases}
\end{equation}
\end{thm}

It is conjectured that the $n$\textsuperscript{th} centered
moments from number theory agree with random matrix theory for any
Schwartz test function; our results above may be interpreted as
providing additional evidence.

Our goal is to reduce as many calculations as possible to ones
already done in the seminal work of \cite{ILS}, where their
delicate analysis of the Kloosterman and Bessel terms in the
Petersson formula allowed them to go well beyond the diagonal. We
quickly review notation and state some needed estimates. We then
calculate the relevant number theory quantities, concentrating on
the new terms that did not arise in \cite{ILS}. Using properties
of the Fourier and Mellin transforms and convolutions, we reduce
our multidimensional integrals of Kloosterman-Bessel terms to
one-dimensional integrals considered in \cite{ILS}.

\begin{rek}
Random matrix theory provides exact formulas for the moments for
test functions of any support, derivable from the $n$-level
densities (in particular, the determinant expansions of these); see
\cite{KS1,KS2} for details. However, \textit{a priori} it is not
obvious that these results agree with those obtained in number
theory for test functions as restricted in our theorems. For
example, much of the analysis in Rubinstein \cite{Ru} and Gao
\cite{Gao} of the $n$-level densities of quadratic Dirichlet
$L$-functions is devoted to analyzing the resulting combinatorial
expressions to show agreement with random matrix theory; the
centered moments are combinatorially much easier to analyze.


To simplify showing agreement between number theory and random
matrix theory we further develop the combinatorics used in the work
of Hughes-Rudnick \cite{HR1} and Soshnikov \cite{Sosh}, and, by
desymmetrizing certain integrals which arise, derive some needed
Fourier transform identities. Doing so allows us to handle support
in $[-\frac1{n-1},\frac1{n-1}]$ on the random matrix side. While
this makes our results more restrictive than the exact determinant
expansions of Katz-Sarnak, these new formulas are significantly more
convenient for comparisons with number theory, involving simple
one-dimensional integrals of convolutions of our test function
rather than sums of determinants. This allows us to avoid the
combinatorial analysis of the number theory terms in \cite{Ru,Gao}.
In the course of proving the agreement between number theory and
random matrix theory, we derive \eqref{eq:thmextmomcompact}, a new
and, for test functions with $\supp(\hphi) \subset (-\frac1{n-1},
\frac1{n-1})$, significantly more tractable expansion for the
$n$\textsuperscript{th} centered moments than the determinant
expansions.\end{rek}

In \S\ref{sect:RMT moments} we see that the first natural boundary
in analyzing the $n$\textsuperscript{th} centered moment for
$\SO({\rm even})$ and $\SO({\rm odd})$ in random matrix theory is
for test functions supported in $[-\frac1n,\frac1n]$; the next
natural boundary (where new terms arise) occurs for test functions
supported in $[-\frac1{n-1},\frac1{n-1}]$. It is essential that we
are able to perform the number theory analysis for test functions
whose support exceeds $[-\frac1n,\frac1n]$. While it is desirable to
obtain as large support as possible, by breaking
$[-\frac1n,\frac1n]$ in Theorem \ref{thm:extended moments} we see
the new terms arise in the number theory expansions as well, and
agree perfectly with random matrix theory.

Instead of investigating centered moments we could study the
$n$-level densities. Assuming GRH, the imaginary parts of the
zeros of an $L$-function associated to a modular form $f \in
H_k^+(N)$  can be written as $\g_f^{(j)}$ where
$0\leq\g_f^{(1)}\leq\g_f^{(2)} \leq\dots$, and $\g_f^{(-j)} =
-\g_f^{(j)}$. If $f\in\hkmn$ there is an additional zero
$\g_f^{(0)}=0$ (note there is no forced zero at the central point
for $f\in\hkpn$). The symmetrized $n$-level density is
\begin{eqnarray}
\frac{1}{|H_k^\pm(N)|} \sum_{f\in H_k^\pm(N)} \sum_{\substack{
j_1,\dots, j_n \\ j_i \neq \pm j_k}} \phi_1\left(\frac{\log
R}{2\pi}\gamma_f^{(j_1)}\right)\cdots \phi_n\left(\frac{\log
R}{2\pi}\gamma_f^{(j_n)}\right),
\end{eqnarray} where the $\phi_i$ are even Schwartz functions
whose Fourier transforms  have compact support. Since our families
are of constant sign, we understand the number of zeros at the
central point (unlike, say, for generic one-parameter families of
$L$-functions of elliptic curves). While our arguments immediately
generalize to the case when the $\phi_i$ are not all equal, we chose
to study the $n$\textsuperscript{th} centered moments to facilitate
comparison with random matrix theory in the range where the
Bessel-Kloosterman terms contribute.

Another application of centered moments is in estimating the order
of vanishing of $L$-functions at the central point. Miller
\cite{Mil0} noticed that as $n$ increases, the $n$-level densities
provide better and better estimates for bounding the order of
vanishing at the central point; unfortunately, as $n$ increases the
bounds are better only for \emph{high} (growing with $n$) vanishing
at the central point. We obtain similar bounds from the $n^{\rm th}$
centered moments. Explicitly, from Theorem \ref{thm:extended
moments} we immediately obtain

\begin{cor}\label{cor:boundsordvancp} Consider the families of weight $k$
cuspidal newforms split by sign, $H_k^\pm(N)$. Assume GRH for all
Dirichlet $L$-functions and all $L(s,f)$. For each $n$ there are
constants $r_n$ and $c_n$ such that as $N\to\infty$ through the
primes, for $r \ge r_n$ the probability of at least $r$ zeros at the
central point is at most $c_n r^{-n}$; equivalently, the probability
of fewer than $r$ zeros at the central point is at least $1-c_n
r^{-n}$.
\end{cor}

The paper is organized as follows. In \S\ref{sec:numthprel} we
review the needed number theory results, write down the expansions
for the centered moment sums, and collect many of the estimates that
we need later. We then prove our number theory results, Theorems
\ref{thm:nosplitmoments} (the unsplit case) and
\ref{thm:mock-Gaussian for D} (the split case with restricted
support) in \S\ref{sec:proofsmockgaussthms}, and Theorem
\ref{thm:extended moments} (the unsplit case where we go beyond the
diagonal by analyzing the Bessel-Kloosterman terms in the Petersson
expansion) in \S\ref{sec:NTncentmom}; we show these agree with
random matrix theory (Theorem \ref{thm:mmts_RMT}) in \S\ref{sect:RMT
moments}. In \S\ref{sec:ordervanLFn} we prove Corollary
\ref{cor:boundsordvancp}, and show that it provides better bounds
than \cite{ILS} for the percentage of odd cuspidal newforms with at
least 5 zeros at the central point.


\section{Number Theory Preliminaries}\label{sec:numthprel}

\subsection{Notation}

\begin{defi}[Gauss Sums]
For $\chi$ a character modulo $q$ and $e(x)  =  e^{2\pi \i x}$,
\begin{equation}\label{eq:gausssum}
G_\chi(n)\ = \ \sum_{a \bmod q} \chi(a) e(an/q),
\end{equation} and $|G_\chi(n)| \le \sqrt{q}$.
\end{defi}

\begin{defi}[Ramanujan Sums]\label{def:RamSums}
If $\chi = \chi_0$ (the principal character modulo $q$) in
\eqref{eq:gausssum}, then $G_{\chi_0}(n)$ becomes the Ramanujan sum
\be\label{eq:defn Rnew} R(n,q) \ = \ \sideset{}{^\ast}\sum_{a \bmod
q} e(an/q) \ = \ \sum_{d|(n,q)} \mu(q/d) d,\ee where $\ast$
restricts the summation to be over all $a$ relatively prime to
$q$.
\end{defi}

\begin{defi}[Kloosterman Sums] For integers $m$ and $n$,
\begin{equation}
S(m,n;q)\ =\ \sideset{}{^\ast}\sum_{d \bmod q} e\left(\frac{md}{q} +
\frac{n\db}{q}\right),
\end{equation}
where $d \db \equiv 1 \bmod q$. We have
\begin{equation}\label{eq:estimate Kloosterman}
|S(m,n;q)|\ \ \leq\ \ (m,n,q)\ \sqrt{\min\left\{ \frac{q}{(m,q)} ,
\frac{q}{(n,q)} \right\}}\ \ \tau(q),
\end{equation}
where $\tau(q)$ is the number of divisors of $q$; see Equation
2.13 of \ils.
\end{defi}

\begin{defi}[Fourier Transform] We use the following
normalization:
\begin{equation}
\hphi(y)\ =\ \intinf \phi(x) e^{-2\pi\i xy} \;\d x, \ \ \ \ \
\phi(x)\ =\ \intinf \hphi(y) e^{2\pi\i xy} \;\d y.
\end{equation}
\end{defi}

\begin{defi}[Characteristic Function]
For $A \subset \R$, let  \begin{equation} \twocase{\I_{\{x \in
A\}} \ = \ }{1}{if $x\in A$}{0}{otherwise.}
\end{equation}
\end{defi}

The Bessel function of the first kind occurs frequently in this
paper, and so we collect here some standard bounds for it (see, for
example, \cite{GR,Wat}).

\begin{lem}\label{lem:Bessel} Let $k\geq 2$ be an
integer. The Bessel function satisfies
\begin{enumerate}
\item\label{lb:1} $J_{k-1}(x) \ll 1$.
\item\label{lb:2} $J_{k-1}(x) \ll x$.
\item\label{lb:3} $J_{k-1}(x) \ll x^{k-1}$.
\item\label{lb:4} $J_{k-1}(x) \ll x^{-\foh}$.
\item\label{lb:5} $2J'_\nu(x) = J_{\nu-1}(x) - J_{\nu+1}(x)$.
\end{enumerate}
\end{lem}

\subsection{Fourier coefficients}

Let $k$ and $N$ be positive integers with $k$ even and $N$ prime.
We denote by $S_k(N)$ the space of all cusp forms of weight $k$
for the Hecke congruence subgroup $\Gamma_0(N)$ of level $N$. That
is, $f$ belongs to $S_k(N)$ if and only if $f$ is holomorphic in
the upper half-plane, satisfies
\begin{equation}
f\left(\frac{az+b}{cz+d}\right)\ =\ (cz+d)^k f(z)
\end{equation}
for all $\left(\smallmatrix a&b\\ c&d\endsmallmatrix\right)\in
\Gamma_0(N) := \left\{\left(\smallmatrix \ga&\gb\\
\gamma&\gd\endsmallmatrix\right) \ : \ \gamma\equiv 0 \bmod
N\right\}$, and vanishes at each cusp of $\Gamma_0(N)$. See \cite{I}
for more details about cusp forms.

Let $f \in S_k(N)$ be a cuspidal newform of weight $k$ and level
$N$; in our case this means $f$ is a cusp form of level $N$ but
not of level $1$. It has a Fourier expansion
\begin{equation}
f(z) \ = \ \sum_{n=1}^\infty a_f(n) e(nz),
\end{equation}
with $f$ normalized so that $a_f(1) = 1$. We  normalize the
coefficients by defining
\begin{equation}
\gl_f(n)\ =\ a_f(n) n^{ -(k-1)/2 }.
\end{equation} $\hkn$ is the set of all
$f\in S_k(N)$ which are newforms. We split this set into two
subsets, $H_k^+(N)$ and $H_k^-(N)$, depending on whether the sign of
the functional equation of the associated $L$-function (see
\S\ref{sect:introduction} for details) is plus or minus. {From}
Equation $2.73$ of \cite{ILS} we have for $N > 1$ that
\begin{equation}\label{eq:number of terms in hkpm}
|\hkpmn| \ = \ \frac{k-1}{24}N + O\left( (kN)^{\frac{5}{6}}
\right).
\end{equation}
For simplicity we shall deal only with the case $N$ prime, a fact
which we will occasionally remind the reader of (though, as in
\cite{ILS}, similar arguments work for $N$ square-free). For a
newform of level $N$, $\gl_f(N)$ is related to the sign of the form
(\cite{ILS}, Equation 3.5):

\begin{lem}\label{lem:sign in terms of lambdaN}
If $f\in\hkn$ and $N$ is prime, then
\begin{equation}\label{eq:signfneqexpils}
\gep_f \ = \ -\i^k \gl_f(N) \sqrt{N}.
\end{equation}
\end{lem}
As $\gep_f = \pm 1$, \eqref{eq:signfneqexpils} implies $|\gl_f(N)| =
\frac1{\sqrt{N}}$. Essential in our investigations will be the
multiplicative properties of the Fourier coefficients.

\begin{lem}\label{lem:multiplicativity of fourier coeffs}
Let $f \in \hkn$. Then
\begin{equation}
\gl_f(m) \gl_f(n) \ = \ \sum_{\substack{d|(m,n) \\ (d,N) = 1}}
\gl_f\left( \frac{mn}{d^2} \right).
\end{equation}
In particular, if $(m,n) = 1$ then
\begin{equation}
\gl_f(m) \gl_f(n)  \ = \  \gl_f(mn),
\end{equation}
and if $p$ is a prime not dividing the level $N$, then
\begin{gather}\label{eq:expansion of lambda_f(p)^n}
\gl_f(p)^{2m} \ = \ \sum_{r=0}^m \left(\binom{2m}{m-r} -
\binom{2m}{m-r-1} \right) \gl_f(p^{2r})\nonumber\\
\gl_f(p)^{2m+1} \ = \ \sum_{r=0}^m \left(\binom{2m+1}{m-r} -
\binom{2m+1}{m-r-1} \right) \gl_f(p^{2r+1}).
\end{gather}
\end{lem}
We discovered the coefficients for the expansion of $\gl_f(p)^n$
from \cite{Guy}. Note for a prime $p\notdiv N$,
\begin{equation}\label{eq:lambda(p) squared}
\gl_f(p)^2 \ = \ \gl_f(p^2) + 1,
\end{equation}
and $\gl_f(p^{2m})$ is a sum of $\gl_f(p^{2r})$ (i.e., only even
powers) while $\gl_f(p^{2m+1})$ is a sum of $\gl_f(p^{2r+1})$ (i.e.,
only odd powers). Consider
\begin{equation}\label{eq:deltakndefeq}
\Delta_{k,N}^\sigma(n) \ = \ \sum_{f \in \hksn} \gl_f(n), \quad
\sigma \in \{+,-,\ast\}.
\end{equation}
Note we are \emph{not} dividing by the cardinality of the family,
which is of order $N$. Splitting by sign and using
Lemma~\ref{lem:sign in terms of lambdaN} we have that if $N$ is
prime and $(N,n)=1$,
\begin{align}
\Delta_{k,N}^\pm(n) \ &= \ \sum_{f\in\hkn} \frac12(1\pm\epsilon_f) \gl_f(n) \nonumber\\
\ &= \ \frac12 \Delta_{k,N}^\ast(n)\ \mp\ \frac{\i^k \sqrt{N}}{2}
\Delta_{k,N}^\ast(nN). \label{eqdeltapm}
\end{align}
Thus, to execute sums over $f\in\hkpmn$, it suffices to understand
sums over all $f\in H_k^\ast(N)$. Propositions 2.1, 2.11 and 2.15
of \cite{ILS} yield a useful form of the Petersson formula:
\begin{lem}[\ils]\label{lemcomplsum}
Let $X, Y$ be parameters to be determined later subject to $X<N$.
If $N$ is prime and $(n,N^2)|N$ then
\begin{equation}\label{eqdeltastarprimeinf}
\Delta_{k,N}^\ast(n)  \ = \  \Delta_{k,N}'(n) +
\Delta_{k,N}^\infty(n) ,
\end{equation}
where
\begin{multline}\label{proptof}
\Delta_{k,N}'(n)\ =\ \frac{(k-1)N}{12\sqrt{n}}\
\delta_{n,\Box_Y} \\
+ \frac{(k-1)N}{12} \sum_{\substack{(m,N) = 1 \\
m \le Y}} \frac{2\pi \i^k}{m} \sum_{\substack{ c\equiv 0 \bmod N
\\ c\geq N}} \frac{ S(m^2,n;c)}{c} \jk{m^2n}{c},
\end{multline}
where $\delta_{n,\Box_Y} = 1$ only if $n=m^2$ with  $m\leq Y$ and
$0$ otherwise. The remaining piece, $\Delta_{k,N}^\infty(n)$, is
called the \emph{complementary sum}.

If $(a_q)$ is a sequence satisfying
\begin{equation}
\sum_{\substack{(q,nN)=1 \\ q< Q}} \gl_f(q) a_q \ \ll\
(nkN)^{\epsilon'}, \ \ \ \log Q \ \ll \ \log kN
\end{equation}
for every\footnote{We need $f\in H_k^\ast(1)$ (as well as $f\in
H_k^\ast(N)$, as these $f$ arise in the combinatorics in expanding
the $\Delta_{k,N}^\ast$.} $f\in H_k^\ast(1) \cup H_k^\ast(N)$, the
implied constant depending on $\gep'$ only, if $(n,N^2)|N$, then by
Lemma 2.12 of \cite{ILS}
\begin{equation}\label{proptot}
\sum_{\substack{(q,nN)=1 \\ q < Q}} \Delta_{k,N}^\infty(nq) a_q \
\ll \ \frac{k N}{\sqrt{(n,N)}} \left( \frac{1}{X} +
\frac{1}{\sqrt{Y}} \right) (nkNXY)^{\gep'}.
\end{equation}
\end{lem}

In the applications we will take $X$ to be either $N-1$ or $N^\gep$
and $Y = N^{\gep}$, where $\gep,\gep'$ are chosen so that the right
hand side of \eqref{proptot} is $O(N^{1-\gep''})$ for some
$\gep''>0$ if $n \notdiv N$, and is $O(N^{-\gep''})$ if $n|N$. In
Lemma \ref{lemcompsumsmall} we show that the complementary sum does
not contribute for all cases that arise in this paper. We write $c =
bN$ for $c \equiv 0 \bmod N$.

Using the estimate on Kloosterman sums, \eqref{eq:estimate
Kloosterman}, the bounds on the Bessel function $J_{k-1}(x) \ll x$
and $J_{k-1}(x) \ll x^{k-1}$ from Lemma~\ref{lem:Bessel}, and
\eqref{eq:number of terms in hkpm}, we can trivially estimate
$\frac{1}{|\hkn|} \Delta_{k,N}'(n)$. We obtain the following
lemma:
\begin{lem}
Assume $(n,N)=1$. Then
\begin{equation}\label{eq:deltaprime}
\frac{1}{|\hkn|}\ \Delta_{k,N}'(n) \ = \ \frac{1}{\sqrt{n}}\
\delta_{n,\Box_Y} + O\left(n^{(k-1)/2} N^{-k+1/2+\gep}\right) ,
\end{equation}
and
\begin{equation}\label{eq:219}
\frac{1}{|\hkn|}\ \Delta_{k,N}'(Nn) \ \ll \ \sqrt{n}
N^{-\frac{3}{2}+\gep} .
\end{equation}
\end{lem}

\begin{rek}\label{rek:whynoharm}
We chose to uniformly average over $f\in H_k^\ast(N)$ in
\eqref{eq:deltakndefeq}. We obtain similar results if instead we use
harmonic averaging as in \cite{Ro} or Theorems 1.9 and 1.10 of
\cite{ILS}, specifically \be \langle Q(f) \rangle_{\pm,{\rm
harmonic}} \ = \ \sum_{f\in H_k^\pm(N)}
\frac{\Gamma(k-1)}{(4\pi)^{k-1} (f,f)_N}\ Q(f), \ee where $(f,f)_N$
is the Petersson inner product on cusp forms of weight $k$ and level
$N$. The advantage of harmonic averaging is that it facilitates the
analysis of the $p^2$ terms in the explicit formula; specifically,
we would not need to assume GRH for $L(s,f)$. We have chosen to use
uniform averages for several reasons. The first is that, as in
Theorem 1.1 of \cite{ILS}, the assumption of GRH for $L(s,f)$ can be
removed relatively easily by appealing to either the Petersson
formula or properties of $L(s, \sym^2 f \otimes \sym^2 f)$. The
second is that much effort was spent in \cite{ILS} in removing these
arithmetic weights (see their comment on page 66), and removing the
weights is essential to bound the order of vanishing at the central
point (see Corollary \ref{cor:boundsordvancp} and Remark
\ref{rek:HoffLocknoharm}). Finally, when we uniformly average, our
transformation of the multidimensional integrals lead to
one-dimensional integrals that are directly comparable to the
uniformly averaged cases in \cite{ILS}.
\end{rek}

The one-dimensional integral referred to above is:
\begin{lem}\label{lem:ilschap7} Let $\Psi$ be an even Schwartz
function with $\supp(\widehat{\Psi}) \subset (-2,2)$. Then
\begin{multline}
\sum_{m \le N^\gep} \frac{1}{m^2} \sum_{(b,N)=1}
\frac{R(m^2,b)R(1,b)}{\varphi(b)}\int_{y=0}^\infty J_{k-1}(y)
\widehat{\Psi}\left( \frac{2 \log (by\sqrt{N} / 4\pi m) }{\log
R}\right) \frac{\d y}{\log R} \\
\ = \ -\frac{1}{2}\left[ \int_{-\infty}^\infty \Psi(x) \frac{\sin
2\pi x}{2\pi x} \d x - \foh \Psi(0)\right] + O\left( \frac{k
\log\log kN}{\log kN}\right),
\end{multline}
where $R(n,c)$ is given by \eqref{eq:defn Rnew}, $R = k^2 N$ and
$\varphi$ is Euler's totient function.
\end{lem}

This follows from Equations 7.5 and 7.6 of \cite{ILS}. The
explicit formula converts sums over zeros to sums over primes.
Later we convert these prime sums to integrals, and then the above
lemma allows us to evaluate the final expressions.

\subsection{Density and Moment Sums}\label{sect:intro_mmt sums}

Let $f\in \hkn$, and let $\Lambda(s,f)$ be its associated
completed $L$-function, \eqref{eq:completed_L_func}. The
Generalized Riemann Hypothesis states that all the zeros of
$\Lambda(s,f)$ (i.e., the non-trivial zeros of $L(s,f)$) are of
the form $\rho_f = \tfrac12 + \i\g_f$ with $\g_f \in \R$. The
analytic conductor of $\Lambda(s,f)$ is $R = k^2 N$, and its
smooth counting function (also called the $1$-level density) is
\begin{eqnarray}\label{eq:expformexpsumzerosdfphi}
D(f;\phi) & \ = \ & \sum_{\gamma_f} \phi\left( \frac{\log R}{2\pi}
\gamma_f \right),
\end{eqnarray}
where $\phi$ an even Schwartz function whose Fourier transform has
compact support and the sum is over all the zeros of $\Lambda(s,f)$.
Because $\phi$ decays rapidly, the main contribution to
\eqref{eq:expformexpsumzerosdfphi} is from zeros near the central
point. The explicit formula applied to $D(f;\phi)$ gives (see
Equation 4.25 of \cite{ILS})
\begin{equation}\label{eqdfphiexpansion}
D(f;\phi)\ =\ \widehat{\phi}(0) + \foh \phi(0) - P(f;\phi) +
O\left( \frac{\log \log R}{\log R} \right),
\end{equation}
where
\begin{equation}\label{eq:P in terms of gl}
P(f;\phi)\ =\ \sum_{p \notdiv N} \gl_f(p) \phir{} \pfrac{}.
\end{equation} While the derivation of \eqref{eqdfphiexpansion} in
\cite{ILS} uses GRH for $L(s,\sym^2 f)$, as they remark this formula
can be established on average over $f$ by an analysis of the
Petersson formula or from properties of $L(s, \sym^2 f \otimes
\sym^2 f)$ (see page 88 of \cite{ILS}). For ease of exposition we
shall assume GRH for $L(s,f)$ below. We trivially absorbed the $p=N$
term into the error. If $\mbox{supp}(\hphi) \subset (-1,1)$,
\cite{ILS} show the $P(f,\phi)$ term does not contribute, and hence
$\lim_{N\to\infty} \<D(f;\phi)\>_\sigma = \hphi(0)+\foh \phi(0)$ for
any $\sigma \in \{+,-,\ast\}$. Thus, to study the centered moments,
we must evaluate
\begin{align}
\left\langle\left(D(f;\phi)-\<D(f;\phi)\>_\sigma\right)^n\right\rangle_\sigma
&\ = \ \left\langle\left(-P(f;\phi) +
\ O\left( \frac{\log \log R}{\log R}\right)\right)^n\right\rangle_\sigma \nonumber\\
&\ = \ (-1)^n\<P(f;\phi)^n\>_\sigma \ + \ O\left( \frac{\log \log
R}{\log R} \right). \label{eq:centered mmts of D in terms of P}
\end{align}
The last line follows from H\"older's inequality and the fact that
$\< P(f;\phi)^n \>_\sigma \ll 1$ (which follows from
\eqref{eqdfphiexpansion} and that $\<|D(f;\phi)|\>_\sigma\ll 1$).
By using H\"older's inequality, we can prove \eqref{eq:centered
mmts of D in terms of P} without having to construct a positive
majorizing test function with suitable support, as is often done
(see, for example, \cite{RS,Ru}). See Appendix
\ref{sec:handerrmomexp} for details. We split by sign and use
Lemma~\ref{lem:sign in terms of lambdaN} to obtain
\begin{align}
\sum_{f\in\hkpmn} P(f;\phi)^n &\ = \
\sum_{f\in\hkn} \frac{1\pm \gep_f}{2} P(f;\phi)^n \nonumber\\
& \ = \ \frac12 \sum_{f\in\hkn}P(f;\phi)^n\ \mp\ \frac{1}{2}
\sum_{f\in\hkn} \i^k \sqrt{N}\gl_f(N)P(f;\phi)^n .
\end{align}
Since $|\hkpn| \sim |\hkmn| \sim \frac12 |\hkn|$ as $N\to\infty$
by \eqref{eq:number of terms in hkpm}, we have
\begin{equation}\label{eq:splitting}
\<P(f;\phi)^n\>_\pm\ \sim\ \< P(f;\phi)^n\>_\ast\ \mp\ \i^k \sqrt{N}
\< \gl_f(N)P(f;\phi)^n \>_\ast.
\end{equation}

In conclusion, if $\supp(\hphi) \subset (-1,1)$, we have
\begin{equation}\label{eq:full_grp_in_terms_of_S1}
\lim_{N\to\infty}
\left\langle\left(D(f;\phi)-\<D(f;\phi)\>_\ast\right)^n\right\rangle_\ast
\ =\ (-1)^n \lim_{N\to\infty} S_1^{(n)}
\end{equation}
and
\begin{equation}\label{eq:nth_mmt_in_terms_of_S}
\lim_{N\to\infty}
\left\langle\left(D(f;\phi)-\<D(f;\phi)\>_\pm\right)^n\right\rangle_\pm
\ =\ (-1)^n \lim_{N\to\infty} S_1^{(n)}\ \pm\ (-1)^{n+1}
\lim_{N\to\infty} S_2^{(n)}
\end{equation}
(assuming all limits exist), where
\begin{equation}\label{eq:S_1}
S_1^{(n)}\ :=\ \sum_{p_1 \notdiv N , \dots, p_n \notdiv N}
\prod_{j=1}^n \left( \phir{j} \pfrac{j} \right) \< \prod_{j=1}^n
\gl_f(p_i) \>_\ast
\end{equation}
and
\begin{equation}\label{eq:S_2}
S_2^{(n)}\ :=\ \i^k \sqrt{N} \sum_{p_1 \notdiv N , \dots, p_n
\notdiv N} \prod_{j=1}^n \left( \phir{j} \pfrac{j} \right) \<
\gl_f(N) \prod_{j=1}^n \gl_f(p_i) \>_\ast.
\end{equation}


\section{Mock-Gaussian behavior:
Proof of Theorems~\ref{thm:nosplitmoments} and
\ref{thm:mock-Gaussian for D}}\label{sec:proofsmockgaussthms}

In this section we prove Theorems~\ref{thm:nosplitmoments} and
\ref{thm:mock-Gaussian for D}, which states that for test
functions with suitable support, the centered moments of
$D(f;\phi)$ are Gaussian. By \eqref{eq:nth_mmt_in_terms_of_S} we
must therefore study the limits of $S_1^{(n)}$ and $S_2^{(n)}$ as
$N\to\infty$ through the primes, with $\supp(\hphi) \subseteq
(-\frac{1}{n} , \frac{1}{n})$.

Because there is no $S_2^{(n)}$ term when we do not split by sign,
Theorem \ref{thm:nosplitmoments} is equivalent to the following
lemma, which we now prove.

\begin{lem}\label{lem:S1}
Let $S_1^{(n)}$ be defined as in \eqref{eq:S_1}, and assume GRH
for $L(s,f)$. Then if $\supp(\hphi) \subset (-\frac{1}{n}
\frac{2k-1}{k} , \frac{1}{n} \frac{2k-1}{k})$,
\begin{equation}
\twocase{\lim_{\substack{N\to\infty \\ N\ {\rm prime}}} S_1^{(n)}
\ = \ }{(2m-1)!!\ \sigma^{2m}_\phi}{if $n=2m$ is even}{0}{if $n$
is odd,}
\end{equation}
where
\begin{equation}
\sigma^2_\phi \ = \ 2\int_{-\infty}^{\infty} |y| \hphi(y)^2 \;\d
y.
\end{equation}
\end{lem}

\begin{proof}
We split the sum over primes into sums over powers of distinct
primes. Let $p_1\cdots p_n = q_1^{n_1} \cdots q_\ell^{n_\ell}$
with the $q_j$ distinct,  so
\begin{equation}
\prod_{j=1}^n \gl_f(p_i)\ =\ \prod_{j=1}^\ell \gl_f(q_j)^{n_j}.
\end{equation}
By the multiplicativity of $\lambda_f$
(Lemma~\ref{lem:multiplicativity of fourier coeffs}),
$\gl_f(q_j)^{n_j}$ can be written as a sum of $\gl_f(q_j^{m_j})$
where the $m_j$ are non-negative integers less than or equals to
$n_j$ with $m_j \equiv n_j \bmod 2$.

The only way for $\prod_{i=1}^n \gl_f(p_i)$ to have a constant
term (i.e., $\gl_f(1)$) is for $p_1\cdots p_n$ to equal a perfect
square; this will be the main term. This can only happen when
$n=2m$ is an even integer. In this case each prime occurs an even
number of times, and the primes can be paired. Assume first that
each $n_j = 2$ so that each prime occurs exactly twice. The number
of ways to pair $2m$ elements in pairs is
$\frac{1}{m!}\binom{2m}{2} \binom{2m-2}{2}\cdots \binom{2}{2} =
\frac{(2m)!}{2^m m!} = (2m-1)!!$; note these are the even moments
of the standard Gaussian. Using the Prime Number Theorem  to
evaluate the prime sums, and the fact that $\hphi$ is even, we see
that the contribution from these terms is
\begin{equation}\label{eq:2mchoosempairs}
\lim_{\substack{N\to\infty \\ N\ {\rm prime}}}(2m-1)!!
\left(\sum_{p\notdiv N} \phir{}^2 \left(\pfrac{}\right)^2
\right)^m\ =\ (2m-1)!! \left(2\int_{-\infty}^\infty \hphi(y)^2 |y|
\; \d y\right)^m;
\end{equation}
note the integral is the variance $\sigma^2_\phi$ because of the
support condition on $\hphi$. The other possibility is that some
$n_j \ge 4$. In this case we obtain a formula similar to
\eqref{eq:2mchoosempairs}, the only changes being a different
combinatorial factor than $(2m-1)!!$ outside, and we have sums such
as \be\label{eq:2mchoosempairsnj} \sum_{p\notdiv N} \phir{}^{n_j}
\left(\pfrac{}\right)^{n_j}. \ee If $n_j = 2$ then
\eqref{eq:2mchoosempairsnj} is $O(1)$ by the Prime Number Theorem;
however, \eqref{eq:2mchoosempairsnj} is $O\left(\log^{-4} R\right)$
if $n_j \ge 4$. Thus the contribution from the terms where at least
one $n_j \ge 4$ is negligible.

The other contributions from expanding $\prod_{i=1}^n \gl_f(p_i)$
are of the form
\begin{equation}
\sum_{\substack{q_1 \notdiv N , \dots, q_\ell \notdiv N \\ q_j
{\rm distinct}}} \prod_{j=1}^\ell \hphi\left(\frac{\log q_j}{\log
R}\right)^{n_j} \left(\frac{2 \log q_j}{\sqrt{q_j} \log R}
\right)^{n_j} \< \gl_f(q_1^{m_1} \cdots q_\ell^{m_\ell}) \>_\ast
\end{equation}
with  $\ell \ge 1$ (i.e., there is at least one prime) and $m_j \geq
1$ for at least one $j$ (i.e., this is not a constant term). We show
in the limit as $N\to\infty$ that these terms do not contribute. By
\eqref{eq:unifaverdef} and \eqref{eqdeltastarprimeinf},
\begin{equation}
\langle \gl_f(q_1^{m_1}\cdots q_\ell^{m_\ell})\rangle_\ast \ = \
\frac{1}{|\hkn|} \left(\Delta_{k,N}'(q_1^{m_1}\cdots
q_\ell^{m_\ell}) + \Delta_{k,N}^\infty(q_1^{m_1}\cdots
q_\ell^{m_\ell}) \right).
\end{equation} Let $X=N-1$ and $Y=N^\gep$.
By Lemma \ref{lemcompsumsmall}, which assumes GRH for $L(s,f)$ (and
in fact is why we assume GRH), for $\gep$ sufficiently small the
complementary sum piece contributes
\begin{equation}
\frac{1}{|\hkn|}\ \Delta_{k,N}^\infty(q_1^{m_1}\cdots
q_\ell^{m_\ell}) \ \ll \ O(N^{-\gep''}).
\end{equation}
For $\Delta'_{k,N}$, by \eqref{eq:deltaprime}
\begin{multline}\label{eq:319s1proof}
\frac{1}{|\hkn|}\ \Delta_{k,N}'(q_1^{m_1}\cdots q_\ell^{m_\ell})\ =\
\frac{1}{q_1^{m_1/2} \cdots q_\ell^{m_\ell/2}}\
\delta_{q_1^{m_1} \cdots q_\ell^{m_\ell},\ \Box_Y} \\
\ + \  O\left((q_1^{m_1} \cdots q_\ell^{m_\ell})^{(k-1)/2}
N^{-k+1/2+\gep}\right),
\end{multline}
where the first term is present only if all $m_j$ are even
(implying all $n_j$ are even as $m_j \equiv n_j \bmod 2$).

First we show the sum over squares is $\ll \log^{-2} R$. The squares
contribute to $S_1^{(n)}$
\begin{equation}
\prod_{j=1}^\ell \sum_{\substack{q_j \\ q_j \neq q_k}}
\hphi\left(\frac{\log q_j}{\log R}\right)^{n_j} \frac{2^{n_j}
\log^{n_j} q_j}{q_j^{(n_j+m_j)/2}\log^{n_j} R}\ .
\end{equation}
Note each $m_j$ is even. The contribution from terms with either
$n_j \ge 2$ and $m_j=0$ or $m_j \ge 2$ is $O(1)$, exactly as above.
However, we have assumed that at least one $m_j\geq 1$ (and since
$m_j$ must be even here, $m_j \geq 2$). The prime sum of such a term
converges, and so its contribution will be $O\left(\log^{-n_j}
R\right)$. The product of all these contributions is at most
$O\left(\log^{-2} R\right)$, as required.

Now we bound the contribution to $S_1^{(n)}$ from the $O$-term in
\eqref{eq:319s1proof}. Recall that $\sum_{j=1}^\ell n_j = n$. If
$\supp(\hphi) \subset (-\ga,\ga)$, the contribution is largest when
$m_j = n_j$, which is bounded by
\begin{align}
& \ll\ N^{-k+\foh+\gep} \sum_{q_1 \ll R^\ga , \dots , q_\ell \ll
R^\ga} \prod_{j=1}^\ell \left[\left(\frac{\log q_j}{\sqrt{q_j}
\log R} \right)^{n_j} q_j^{m_j(k-1)/2} \right] \nonumber\\
& \ll\ N^{-k+\foh+\gep} \prod_{j=1}^\ell \left(\sum_{q \ll R^\ga}
q^{-\frac{n_j k}2-n_j}\right)\nonumber\\
& \ll\ N^{-k+\foh+\gep} R^{\frac{nk\alpha}2-n+\ell}.
\end{align}
The worst case is when $\ell = n$. Since $R=k^2N$, if
$\alpha<\frac{2k-1}{nk}$ this vanishes as $N\to\infty$. Hence if
$\ga < \frac{1}{n}\frac{2k-1}{k}$
\begin{equation}
\lim_{\substack{N\to\infty \\ N\ {\rm prime}}} S_1^{(n)}\ = \
\begin{cases}
(2m-1)!! \left(2\int_{-\infty}^\infty \hphi(y)^2 |y| \; \d
y\right)^m &
\text{ if $n=2m$ is even}\\
0 & \text { if $n$ is odd.}
\end{cases}
\end{equation}
This completes the proof of Lemma \ref{lem:S1} and Theorem
\ref{thm:nosplitmoments}.
\end{proof}

In Theorem \ref{thm:nosplitbetterk}, by assuming GRH for Dirichlet
$L$-functions, we extend the support in Theorem
\ref{thm:nosplitmoments} to $(-\frac2n, \frac2n)$ for $2k \ge n$.

Theorem \ref{thm:mock-Gaussian for D} is equivalent to showing that
$S_2^{(n)}$ is negligible for $\supp(\hphi) \subset (-\frac1n,
\frac1n)$, which we now prove.

\begin{lem}\label{lem:S2}
Assume GRH for $L(s,f)$, and let $S_2^{(n)}$ be defined as in
\eqref{eq:S_2}. If $\supp(\hphi) \subset (-\frac{1}{n} ,
\frac{1}{n} )$, then
\begin{equation}
\lim_{\substack{N\to\infty \\ N\ {\rm prime}}} S_2^{(n)} = 0 .
\end{equation}
\end{lem}

\begin{proof}
The same argument for $S_1^{(n)}$ works for $S_2^{(n)}$, but now
there can be no squares because we have $\gl_f(N)$ and none of the
primes equal $N$. $S_2^{(n)}$ is made up of terms like
\begin{equation}
\sqrt{N} \sum_{\substack{q_1 \notdiv N , \dots, q_\ell \notdiv N \\
q_j {\rm distinct}}}\ \prod_{j=1}^\ell \hphi\left(\frac{\log
q_j}{\log R}\right)^{n_j} \left(\frac{2 \log q_j}{\sqrt{q_j} \log
R} \right)^{n_j} \< \gl_f(N q_1^{m_1} \dots q_\ell^{m_\ell})
\>_\ast.
\end{equation}
We again use Lemma~\ref{lemcomplsum} to evaluate the average over
$\gl_f$. By Lemma~\ref{lemcompsumsmall} (which requires GRH for
$L(s,f)$) the complementary sum is $O(N^{-1-\gep''})$, which is
negligible when multiplied by $N^{1/2}$. By \eqref{eq:219} the
remaining piece is bounded by
\begin{align}
\ll N^{-1+\gep} \sum_{q_1 \ll R^\ga , \dots , q_\ell \ll R^\ga}
\prod_{j=1}^\ell \left[\left(\frac{2\log q_j}{\sqrt{q_j} \log R}
\right)^{n_j} q_j^{m_j/2} \right] &\ll\ N^{-1+\gep} \left(\sum_{q
\ll R^\ga} \frac{\log q}{\log R}\right)^n \nonumber\\
& \ll\ N^{-1+\gep} R^{n\ga},
\end{align}
as the worst term occurs when $n_j = m_j = 1$. This contribution is
vanishingly small if $\ga < \frac{1}{n}$ (recall $R = k^2N$).
\end{proof}

Therefore, by \eqref{eq:full_grp_in_terms_of_S1} and
Lemma~\ref{lem:S1}, if $\supp(\hphi) \subset
(-\frac{1}{n}\frac{2k-1}{k} ,\frac{1}{n}\frac{2k-1}{k})$ then
\begin{equation} \label{eq:mock Gaussian for full group}
\lim_{\substack{N\to\infty \\ N\ {\rm prime}}}
\left\langle\left(D(f;\phi)-\<D(f;\phi)\>_\ast\right)^n\right\rangle_\ast
=
\begin{cases}
(2m-1)!! \left(2\intinf \hphi(y)^2 |y| \; \d y\right)^m & \text{ if $n=2m$ is even}\\
0 & \text { if $n$ is odd,}
\end{cases}
\end{equation}
and by \eqref{eq:nth_mmt_in_terms_of_S} and Lemmas~\ref{lem:S1}
and \ref{lem:S2}, if $\supp(\hphi) \subset (-\frac{1}{n} ,
\frac{1}{n})$ then
\begin{equation}\label{eq:mock Gaussian for split group}
\lim_{\substack{N\to\infty \\ N\ {\rm prime}}}
\left\langle\left(D(f;\phi)-\<D(f;\phi)\>_\pm\right)^n\right\rangle_\pm
=
\begin{cases}
(2m-1)!! \left(2\intinf \hphi(y)^2 |y| \; \d y\right)^m & \text{ if $n=2m$ is even}\\
0 & \text { if $n$ is odd.}
\end{cases}
\end{equation}
Because of the support condition on $\hphi$, the integral in
\eqref{eq:mock Gaussian for split group} is the same as the
variance in Theorem~\ref{thm:mock-Gaussian for D}, which completes
the proof of that theorem.

\begin{rek}\label{rek:suppineqMGfFG}
By choosing $k$ sufficiently large, we can take the support of
$\hphi$ as close to $(-\frac2{n},\frac2{n})$ as desired in Theorem
\ref{thm:nosplitmoments}; by using GRH for Dirichlet $L$-functions
in Theorem \ref{thm:nosplitbetterk} we show that if $k$ is
sufficiently large relative to $n$ ($2k \ge n$) then we may take any
$\hphi$ with $\supp(\hphi) \subset (-\frac2{n},\frac2{n})$. This is
a natural boundary to expect, as \cite{ILS} obtained $(-2,2)$ when
$n=1$. For mock-Gaussian behavior (Theorem~\ref{thm:mock-Gaussian
for D}), we do not need to be able to handle support as large as
that; however, support exceeding $(-\frac1n,\frac1n)$ will be
essential in calculating the centered moments in the extended regime
of Theorem~\ref{thm:extended moments}.
\end{rek}


\section{Going Beyond the Diagonal:
Proof of Theorem~\ref{thm:extended moments}}\label{sec:NTncentmom}

We calculate the $n$\textsuperscript{th} centered moment of
$D(f;\phi)$ when $n\geq 2$ and $\supp(\hphi) \subset (-\frac1{n-1},
\frac1{n-1})$, and we will not worry about terms which do not
contribute in this region.  We outline the arguments below. We
assume GRH for $L(s,f)$ for ease of exposition, though as stated in
Remark \ref{rek:whyGRHforLsf}, following \cite{ILS} we may remove
this assumption with additional effort. In \S\ref{sec:14prelims} we
reduce the proof of Theorem \ref{thm:extended moments} to the limit
of $S_2^{(n)}$, which we analyze in the following subsections. In
\S\ref{sec:14peters} we apply the Petersson formula, and in
\S\ref{sec:kloostermansumexpansion} we analyze the Kloosterman terms
by using Dirichlet characters. In
\S\ref{sec:handlingnontrivialchars} we see that the contributions
from the non-principal characters are negligible. By using the
Mellin transform and shifting contours, we convert the prime sums to
integrals in Lemma \ref{lem:doing_the_prime_sums} in
\S\ref{sec:14convert}. The proof of Theorem \ref{thm:extended
moments} is completed by evaluating these integrals in
\S\ref{sec:14changevar}, where by changing variables Lemma
\ref{lem:ilschap7} is applicable.

\subsection{Preliminaries}\label{sec:14prelims}
As \cite{ILS} has already handled the case when $n=1$, we assume
$n\ge 2$ below. Let $\supp(\hphi) \subset (-\sigma,\sigma)$ with
$\sigma\leq 1$. By \eqref{eq:nth_mmt_in_terms_of_S},
\begin{equation}
\lim_{N\to\infty}
\left\langle\left(D(f;\phi)-\<D(f;\phi)\>_\pm\right)^n\right\rangle_\pm
\ =\ (-1)^n \lim_{N\to\infty} S_1^{(n)} \pm (-1)^{n+1}
\lim_{N\to\infty} S_2^{(n)}.
\end{equation}
To prove Theorem~\ref{thm:extended moments} we need to handle
support up to $\frac{1}{n-1}$. If $n \ge 3$ and $k\ge 2$, then
$\frac{1}{n-1} \leq \frac{1}{n}\frac{2k-1}{k}$, and thus
Lemma~\ref{lem:S1} evaluates $S_1^{(n)}$ for $\sigma < \frac1{n-1}$.
If $n=2$, however, then $\frac1{n-1} > \frac{1}{n}\frac{2k-1}{k}$,
and thus there is a decrease in support. This is easily surmounted
by using Theorem \ref{thm:nosplitbetterk} instead of Lemma
\ref{lem:S1}. Theorem \ref{thm:nosplitbetterk} assumes GRH for
Dirichlet $L$-functions; however, we shall be assuming GRH for
Dirichlet $L$-functions when we study $S_2^{(n)}$.

Thus all that remains to prove Theorem \ref{thm:extended moments} is
to show that if $\sigma<\frac{1}{n-1}$ then
\begin{equation}
\lim_{\substack{N\to\infty \\ N\ {\rm prime}}} S_2^{(n)}\ =\
2^{n-1}\left[ \int_{-\infty}^\infty \phi(x)^{n} \frac{\sin 2\pi
x}{2\pi x} \;\d x - \foh \phi(0)^{n} \right],
\end{equation}
and this we shall proceed to do in a series of lemmas, culminating
in Lemma~\ref{lem:final answer for S_2}. This will complete the
proof of the $n$\textsuperscript{th} centered moment in
Theorem~\ref{thm:extended moments}.

\begin{rek} When we do not split by sign (as in Theorems
\ref{thm:nosplitmoments} and \ref{thm:nosplitbetterk}), we can prove
results up to $\frac2{n}$; because of the more complicated terms in
the Bessel-Kloosterman expansion, we can only handle the split cases
up to $\frac1{n-1}$. As the two supports are equal when $n=2$,
investigating small $n$ can be quite misleading as to what support
one should expect for general $n$.\end{rek}

\subsection{Applying the Petersson Formula}\label{sec:14peters}

\begin{lem}\label{lem:natural boundary}
Let $S_2^{(n)}$ be defined as in \eqref{eq:S_2}, and assume GRH
for $L(s,f)$. If $\supp(\hphi) \subset (-\frac{1}{n-1} ,
\frac{1}{n-1})$, then
\begin{multline}\label{eq:centeredmomentkeyexpn}
S_2^{(n)}\ =\ \frac{2^{n+1}\pi}{\sqrt{N}} \sum_{p_1, \dots , p_n}
\sum_{m \le N^\epsilon} \frac{1}{m} \sum_{b=1}^\infty
\frac{S(m^2,p_1 \cdots p_n N;Nb)}{b}
J_{k-1}\left(\frac{4\pi m \sqrt{p_1 \cdots p_n}}{b\sqrt{N}}\right)\\
\times  \prod_{j=1}^n \left(\hphi\left(\frac{\log p_j}{\log
R}\right) \frac{\log p_j}{\sqrt{p_j}\log R} \right) +
O(N^{-\epsilon}).
\end{multline}
\end{lem}

\begin{proof}
The multiplicativity of $\lambda_f$
(Lemma~\ref{lem:multiplicativity of fourier coeffs}) shows that
$S_2^{(n)}$ is made up of terms of the form
\begin{equation}\label{eq:S_2(n)_expanded_out}
\i^k \sqrt{N} \sideset{}{'}\sum_{\substack{q_1 \notdiv N , \dots, q_\ell \notdiv N \\
q_j {\rm distinct}}}\ \prod_{j=1}^\ell \hphi\left(\frac{\log
q_j}{\log R}\right)^{n_j} \left(\frac{2 \log q_j}{\sqrt{q_j} \log
R} \right)^{n_j} \< \gl_f(N q_1^{m_1} \cdots q_\ell^{m_\ell})
\>_\ast,
\end{equation}
where $m_j \leq n_j$, $m_j \equiv n_j \mod 2$ and $\sum n_j = n$;
here and below \tiny$\sideset{}{'}\sum$\normalsize\ means the sum is
taken over distinct primes only. We will show that the contribution
from terms with at least one $n_j\geq 2$ is vanishingly small as
$N\to\infty$ when $\supp(\hphi) \subset (-\frac{1}{n-1} ,
\frac{1}{n-1})$.

We expand $\< \gl_f(N q_1^{m_1} \cdots q_\ell^{m_\ell}) \>_\ast$ via
the Petersson formula (Lemma~\ref{lemcomplsum}). By
Lemma~\ref{lemcompsumsmall}, which relies on GRH for $L(s,f)$, the
complementary sums are of size $O(N^{-1-\gep''})$ for $X=Y =
N^\gep$, which is negligible when multiplied by $N^{1/2}$. We are
left with the $\Delta_{k,N}'(N q_1^{m_1} \cdots q_\ell^{m_\ell})$
terms. That is, \eqref{eq:S_2(n)_expanded_out} can be replaced by
\begin{equation}\label{eq:E_def}
E :=  \i^k \sqrt{N} \sideset{}{'}\sum_{\substack{q_1 \notdiv N , \dots, q_\ell \notdiv N \\
q_j {\rm distinct}}}\left[ \prod_{j=1}^\ell \hphi\left(\frac{\log
q_j}{\log R}\right)^{n_j} \left(\frac{2 \log q_j}{\sqrt{q_j} \log
R} \right)^{n_j}\right] \frac{1}{|\hkn|} \Delta'_{k,N}(N q_1^{m_1}
\cdots q_\ell^{m_\ell}) .
\end{equation}

Assume that $\supp(\hphi) \subseteq [-\sigma,\sigma]$. Note that
$N q_1^{m_1} \cdots q_\ell^{m_\ell}$ can never equal a square,
since none of the $q_j$ divide $N$.  Applying \eqref{eq:219} we
obtain
\begin{equation}\label{eq:219_with_N}
\frac{1}{|\hkn|}\ \Delta'_{k,N}(N q_1^{m_1} \cdots q_\ell^{m_\ell})
\ \ll\ N^{-3/2+\epsilon} q_1^{m_1/2} \cdots q_\ell^{m_\ell/2},
\end{equation}
and so
\begin{equation}\label{eq:qwertyuiop}
E\ \ll\ \sum_{q_1 \ll R^\sigma , \dots , q_{\ell} \ll R^\sigma}
\left[\prod_{j=1}^\ell \left(\frac{\log q_j}{\sqrt{q_j} \log R}
\right)^{n_j}\right] \frac{1}{N^{1-\epsilon}} q_1^{m_1/2} \cdots
q_\ell^{m_\ell/2}.
\end{equation}
The sum in \eqref{eq:qwertyuiop} is maximized if $m_j = n_j$ and as
many as possible of the $n_j=1$, because this maximizes $\ell$ and
hence the number of sums. For the cases where at least one $n_j \ge
2$, the worst case is when $\ell=n-1$, whence the sum in
\eqref{eq:qwertyuiop} contributes
\begin{equation}
\frac{1}{( \log R)^n N^{1-\epsilon}} \sum_{q_1 \ll R^\sigma ,
\dots , q_{n-1} \ll R^\sigma } \left(\log q_1 \right) \cdots
\left(\log q_{n-2} \right) \left(\log q_{n-1} \right)^2 \ll
N^{-1+\epsilon} R^{(n-1)\sigma}.
\end{equation}
If $\sigma<\frac{1}{n-1}$ this has a negligible contribution in
the large $N$ limit. Therefore if $\sigma<\frac{1}{n-1}$ the only
way for \eqref{eq:E_def} not to vanish as $N\to\infty$ is if all
the $n_j=1$. In other words we have shown that
\begin{equation}
S_2^{(n)}\ =\ \frac{\i^k \sqrt{N}}{|\hkn|}
\sideset{}{'}\sum_{\substack{p_1 \notdiv N , \dots, p_n \notdiv N\\
p_j \text{ distinct}}} \prod_{j=1}^n \left( \phir{j} \pfrac{j}
\right) \Delta'_{k,N}(N p_1 \cdots p_n) + O(N^{-\epsilon}).
\end{equation}
We remove the distinctness condition by trivially summing the
contribution when two or more primes coincide. If $p_{n-1}=p_n$,
say, then by \eqref{eq:219_with_N} this contributes
\begin{equation}
\ll\ \sum_{p_1, \dots , p_{n-1} \ll R^\sigma}
\left[\prod_{j=1}^{n-2} \left(\frac{\log p_j}{\sqrt{p_j} \log R}
\right)\right] \left(\frac{\log p_{n-1}}{\sqrt{p_{n-1}} \log R}
\right)^2 \frac{1}{N^{1-\epsilon}}\ p_1^{1/2} \cdots p_{n-2}^{1/2}\
p_{n-1},
\end{equation}
which is of size $N^{-1+\epsilon} R^{(n-1)\sigma}$ and is
vanishingly small if $\sigma<1/(n-1)$. Since $R=k^2N$ and $N$ is a
prime, the compact support condition on $\hphi$ means the condition
$p_j \nmid N$ is automatically satisfied for sufficiently large $N$.
Finally, since \eqref{eq:number of terms in hkpm} shows that $|\hkn|
\sim N(k-1)/12$, applying \eqref{proptof} with $X=Y=N^\epsilon$
yields the lemma.
\end{proof}

\begin{rek}
If $\sigma>\frac{1}{n-1}$, the contribution to the
$n$\textsuperscript{th} centered moment arising from powers of
primes needs to be considered; however, other calculations below
(Lemma \ref{lem:doing_the_prime_sums}) can only be analyzed for
$\sigma < \frac1{n-1}$. In \S\ref{sect:RMT moments} we see this is a
natural boundary, and that new terms are expected to arise once the
support exceeds $[-\frac1{n-1},\frac1{n-1}]$.
\end{rek}

\begin{lem}\label{lem:truncate b sum to coprime}
For $\supp(\hphi) \subseteq (-\frac{5}{2n},\frac{5}{2n})$, the
contribution in \eqref{eq:centeredmomentkeyexpn} from the terms
when $(b,N) \neq 1$ is $O(N^{-\gep})$.
\end{lem}

\begin{proof}
Since $N$ is prime, if $(b,N) \neq 1$ then $(b,N) = jN$ for
$j=1,2,\dots$. If $\supp \hphi \subset (-\sigma,\sigma)$ then
these terms contribute to $S_2^{(n)}$ an amount bounded by
\begin{equation}
\ll \frac{1}{\sqrt{N}} \sum_{p_1, \dots , p_n \leq N^\sigma} \sum_{m
\le N^\epsilon} \frac{1}{m} \sum_{j=1}^\infty \frac{|S(m^2,p_1
\cdots p_n N;j N^2)|}{j N} \left|J_{k-1}\left(\frac{4\pi m \sqrt{p_1
\cdots p_n}}{j N^{3/2}}\right)\right| \frac{1}{\sqrt{p_1 \cdots
p_n}}.
\end{equation}
By the bound for Kloosterman sums \eqref{eq:estimate Kloosterman},
$|S(m^2,p_1\cdots p_n N;j N^2)| \ll j^{\foh+\gep} N^{\foh +\gep}$.
This is because $(m^2,p_1\cdots p_nN,j N^2) \ll m^2 \ll N^\gep$
and $\tau(c) \ll c^\gep$. Lemma \ref{lem:Bessel}(\ref{lb:2}) gives
$J_{k-1}(x) \ll x$, and thus the contribution is bounded by
\begin{equation}
\ll N^{-5/2+\gep'} N^{n\sigma} \left(\sum_{m \le N^\epsilon} 1
\right) \left(\sum_{j=1}^\infty  \frac{1}{j^{3/2-\gep}} \right) \ll
N^{-5/2 + n\sigma +\gep''},
\end{equation}
which is vanishingly small if $\sigma<5/2n$.
\end{proof}

Combining Lemmas \ref{lem:natural boundary} and \ref{lem:truncate b
sum to coprime}, we have under GRH for $L(s,f)$, if $\supp(\hphi)
\subset (-\frac1{n-1},\frac{1}{n-1})$ then
\begin{multline}\label{eq:multilineqwertewqn}
S_2^{(n)} = \frac{ 2^{n+1}\pi}{\sqrt{N}} \sum_{p_1, \dots , p_n}
\sum_{m \le N^\epsilon} \frac{1}{m} \sum_{(b,N)=1} \frac{S(m^2,p_1
\cdots p_n N;Nb)}{b}
J_{k-1}\left(\frac{4\pi m \sqrt{p_1 \cdots p_n}}{b\sqrt{N}}\right)\\
\times  \prod_{j=1}^n \left(\hphi\left(\frac{\log p_j}{\log
R}\right) \frac{\log p_j}{\sqrt{p_j}\log R} \right) +
O(N^{-\epsilon}).
\end{multline}

We now show that the terms with $\log b \gg \log N$ are negligible.
We need to restrict $b$ because later (equation \eqref{eq:defnX2d})
we have sums of $1/b$, and this ensures that sum is not too large.

\begin{lem}\label{lem:truncate b sum}
For $\supp(\hphi) \subseteq (-\frac{1000}{n},\frac{1000}{n})$, the
contribution in \eqref{eq:multilineqwertewqn} from the $b \ge
N^{2006}$ terms is $O(N^{-3})$.
\end{lem}

\begin{proof}
By the bound for Kloosterman sums \eqref{eq:estimate Kloosterman},
$S(m^2,p_1\cdots p_nN,bN) \ll b^{\foh+\gep} N^\gep$. This is because
$(m^2,p_1\cdots p_nN,bN) \ll m^2 \ll N^\gep$ and $\tau(c) \ll
c^\gep$. Lemma \ref{lem:Bessel}(\ref{lb:3}) gives $J_{k-1}(x) \ll
x^{k-1}$, which bounds the summand in
\eqref{eq:centeredmomentkeyexpn} by
\begin{multline}
\frac{1}{\sqrt{N}} \frac{1}{m} \frac{b^{\foh+\gep}N^\gep}{b} m^{k-1}
(p_1\cdots p_n)^{\frac{k-1}{2}} N^{-\frac{k-1}{2}}
b^{-(k-1)} \frac{1}{\sqrt{p_1\cdots p_n}} \\
=  m^{k-2} b^{- k+\foh} (p_1\cdots p_n)^{\frac{k}{2} - 1}
N^{-\frac{k}{2}} N^\gep.
\end{multline}
If $\supp(\hphi) \subset [-\sigma,\sigma]$ then the $\phir{j}$ term
in \eqref{eq:centeredmomentkeyexpn} restricts the $p_j$-sum to be
over $p_j \ll N^\sigma$. Executing the summations over the primes
and summing over $b \ge N^{2006}$ yields \bea N^{-\frac{k}{2}+\gep}
\sum_{m \le N^\gep} m^{k-2} \sum_{b \ge N^{2006}} b^{- k+\foh}
\prod_{j=1}^n \sum_{p_j \ll N^\sigma}p_j^{\frac{k}{2} - 1}
 & \ \ll \ & N^{-\frac{k}{2}+\gep'} N^{(-k + \frac{3}{2})2006} N^{\frac{k}{2}n\sigma}
 \nonumber\\ & \ \ll \ & N^{\frac{k}2(n\sigma - 1004 + \gep')}. \eea
 Therefore the contribution as $N\to\infty$ in
\eqref{eq:multilineqwertewqn} from the terms when $b\geq N^{2006}$
is negligibly small when $\sigma<\frac{1000}{n}$.
\end{proof}

\subsection{Expanding the Kloosterman
Sums}\label{sec:kloostermansumexpansion}

The next lemma converts the Kloosterman sums into Gauss sums, and in
\S\ref{sec:14convert} in Lemma \ref{lem:doing_the_prime_sums} we
convert the resulting prime sums into an integral, which we evaluate
in \S\ref{sec:14changevar}, completing the proof of Theorem
\ref{thm:extended moments}.

\begin{lem}\label{lem:E1contrmultiqwrt}
Under GRH for $L(s,f)$, if $\supp(\hphi) \subset
(-\frac1{n-1},\frac{1}{n-1})$ then
\begin{multline}\label{eq:multilineqwertewqn2}
S_2^{(n)} = -\frac{2^{n+1}\pi}{\sqrt{N}} \sum_{p_1, \dots , p_n}
\sum_{m \le N^\epsilon} \frac{1}{m} \sum_{\substack{(b,N)=1 \\ b <
N^{2006}}} \frac{1}{b\varphi(b)} \sum_{\chi \smod{b}} \chi(N)
G_{\chi}(m^2) G_{\chi}(1)
\overline{\chi}(p_1\cdots p_n)\\
\times J_{k-1}\left(\frac{4\pi m \sqrt{p_1 \cdots
p_n}}{b\sqrt{N}}\right) \prod_{j=1}^n \left(\hphi\left(\frac{\log
p_j}{\log R}\right) \frac{\log p_j}{\sqrt{p_j}\log R} \right) +
O(N^{-\gep}).
\end{multline}
\end{lem}

\begin{proof}
By Lemma \ref{lem:usingsmallermodsimp} we have for $(p_1\cdots
p_n,b) = 1$ and $(b,N) = 1$ that
\begin{equation}\label{eq:fromlemmausingsmaller}
S(m^2,p_1\cdots p_n N; Nb) \ = \  \frac{-1}{\varphi(b)} \sum_{\chi
\smod{b}} \chi(N) G_{\chi}(m^2) G_{\chi}(1)
\overline{\chi}(p_1\cdots p_n).
\end{equation}
If $(p_1\cdots p_n,b)> 1$ then the left hand side of
\eqref{eq:fromlemmausingsmaller} is non-zero but the right hand side
vanishes; however, the contribution to \eqref{eq:multilineqwertewqn}
when $(p_1\cdots p_n,b) > 1$ is negligible if $\supp(\hphi) \subset
(-\frac1{n-1},\frac1{n-1})$. To see this, note that the worst case
is when just one prime divides $b$ and the other $n-1$ primes range
freely. We may assume $p_1|b$, and write $b=r p_1$ (since $p_1$ is a
prime). As $J_{k-1}(x) \ll x$, such terms contribute to
\eqref{eq:multilineqwertewqn} an amount bounded by
\begin{align}
\frac{1}{\sqrt{N}} \sum_{p_1,\dots,p_n \leq N^\sigma} \sum_{m\leq N^\gep}
\frac{1}{m} \sum_{r=1}^\infty \frac{m^2 \sqrt{r}}{r p_1}
\frac{\sqrt{p_1\dots p_n}}{r p_1 \sqrt{N}} \frac{1}{\sqrt{p_1 \dots p_n}}
&\ \ll\ N^{-1 + \gep'}   \left(\sum_{p \leq N^\sigma} 1 \right)^{n-1} \nonumber\\
&\ \ll\ N^{-1+(n-1)\sigma + \gep'},
\end{align}
which is vanishingly small if $\sigma<\frac{1}{n-1}$.

Thus we may use \eqref{eq:fromlemmausingsmaller} in
\eqref{eq:multilineqwertewqn} for all $(p_1,\dots,p_n,b,m)$, which
yields the lemma. Note that the minus sign comes from the $-1$ in
\eqref{eq:fromlemmausingsmaller} from Lemma
\ref{lem:usingsmallermodsimp}.
\end{proof}

\subsection{Handling the Non-Principal
Characters}\label{sec:handlingnontrivialchars}

\begin{lem}\label{lem:E1contrmultiqwrt2}
Under GRH for Dirichlet $L$-functions, if $\supp(\hphi) \subseteq
\left(-\frac2{n},\frac2{n}\right)$ then the contribution from the
non-principal characters to $S_2^{(n)}$ in
\eqref{eq:multilineqwertewqn2} is negligible.
\end{lem}

\begin{proof}  We use $J_{k-1}(x) \ll x$ to bound
the contribution from the non-principal characters in
\eqref{eq:multilineqwertewqn2} by
\begin{eqnarray}\label{eq:defnX2c} & \ \ll \ & \frac1{\sqrt{N}} \sum_{m \le N^\gep}
\frac{1}{m} \sum_{\substack{(b,N) = 1 \\ b < N^{2006}}} \frac{1}{b}
\frac1{\varphi(b)}
 \sum_{\smallsubstack{\chi \smod{b}}{\chi \neq \chi_0}}
\left|G_{\chi}(m^2) G_{\chi}(1)\right| \nonumber\\ & & \ \ \ \ \ \ \
\ \times \frac{m}{b\sqrt{N}} \prod_{j=1}^n \left|\sum_{p_j \neq N}
\overline{\chi}(p_j) \log p_{j} \cdot \frac{1}{\log R}
\phir{i}\right|.
\end{eqnarray} As $\chi \neq \chi_0$ (the principal character with
modulus $b$), by GRH for Dirichlet $L$-functions we have for $\log
xNb \ll R$ that $\sum_{p \le x} \overline{\chi}(p) \log p =
O(x^{\foh} \log^2(bNx))$ $=$ $O(x^\foh R^\gep)$. We now use partial
summation and the compact support of $\hphi$. The boundary term
vanishes, and we are left with \bea\label{eq:418touselater}
\sum_{p_j \neq N} \overline{\chi}(p_j) \log p_j \cdot \frac{1}{\log
R} \phir{j} &\ \ll \ & R^\gep \int_2^{R^\sigma} u^\foh
\left|\frac{\d}{\;\d
u}\hphi\left( \frac{\log u}{\log R}\right) \right| \;\d u \nonumber\\
&\ \ll \ & R^\gep \int_2^{R^\sigma} u^{-\foh} \left|\frac1{\log
R}\hphi'\left( \frac{\log u}{\log R}\right)\right|
 \;\d u\nonumber\\ & \ \ll \ & R^{(\frac{1}2+\gep)\sigma}. \eea As $R =
 k^2N$, the contribution from the $n$ prime sums in \eqref{eq:multilineqwertewqn2}
 is $\ll N^{\frac{\sigma n}{2}+ \gep'}$.

By Lemma \ref{lem:csboundchar}, \be \frac1{\varphi(b)}
 \sum_{\smallsubstack{\chi \smod{b}}{\chi \neq \chi_0}}
\left|G_{\chi}(m^2) G_{\chi}(1)\right| \ \ll \ b. \ee Substituting
the character and prime sum bounds into \eqref{eq:defnX2c} and
executing the sum on $m$ yields \be\label{eq:defnX2d} N^{-\foh}
N^{\gep''} \sum_{\substack{(b,N) = 1 \\ b < N^{2006}}}
\frac{b}{b^2\sqrt{N}}\ N^{\frac{\sigma n}2} \ \ll \ N^{-1 +
\frac{n}{2}\sigma + \gep'''}. \ee As we have $\sum b^{-1}$ above, it
is essential that $b$ is at most a fixed power of $N$; this is
accomplished by Lemma \ref{lem:truncate b sum}. Therefore the
non-principal characters do not contribute to
\eqref{eq:multilineqwertewqn2} for $\supp(\hphi) \subset
\left(-\frac2{n},\frac2{n}\right)$.
\end{proof}

The next lemma shows if
$\supp(\hphi)\subset(-\frac{1}{n-1},\frac{1}{n-1})$, then we may add
in the contribution of powers of primes with negligible error. This
aids the passage to $L$-functions in the next section.

\begin{lem}\label{lem:expand_S_2_with_prim_char}
Under GRH for $L(s,f)$ and all Dirichlet $L$-functions, if
$\supp(\hphi) \subset (-\frac1{n-1},\frac1{n-1})$, then
\begin{multline}\label{eq:multilineqwertewqn5}
S_2^{(n)} = - \frac{2^{n+1}\pi}{\sqrt{N}} \sum_{m \le N^\epsilon}
\frac{1}{m} \sum_{\substack{(b,N) = 1 \\ b < N^{2006}}}
\frac{R(m^2,b)R(1,b)}{b\varphi(b)} \\
\times \sum_{n_1,\dots,n_n} \left(\prod_{j=1}^n
\hphi\left(\frac{\log n_j}{\log R}\right)
\frac{\chi_0(n_j)\Lambda(n_j)}{\sqrt{n_j}\log R} \right)
J_{k-1}\left(\frac{4\pi m \sqrt{n_1\dots n_n}}{b\sqrt{N}}\right) +
O(N^{-\epsilon}).
\end{multline}
\end{lem}

\begin{proof} The support condition follows from taking the minimum
of the supports of Lemmas \ref{lem:natural boundary},
\ref{lem:truncate b sum to coprime}, \ref{lem:truncate b sum},
\ref{lem:E1contrmultiqwrt} and \ref{lem:E1contrmultiqwrt2}, and
\eqref{eq:addsqhigherpowers} below.

Let $\chi_0$ be the principal character modulo $b$. Since it is
real, $\overline{\chi_0} = \chi_0$. {From} \eqref{eq:defn Rnew},
the definition of $R(\alpha,b)$, we have
$R(\alpha,b)=G_{\chi_0}(\alpha)$. Thus $G_{\chi_0}(m^2)
G_{\chi_0}(1) = R(m^2,b) R(1,b)$. Since $(b,N)=1$, $\chi_0(N) =
1$. Lemmas \ref{lem:E1contrmultiqwrt} and
\ref{lem:E1contrmultiqwrt2} imply \eqref{eq:multilineqwertewqn5},
with the restriction that the sums are taken over primes.

We must show that the squares and higher powers of the primes add a
negligible contribution to \eqref{eq:multilineqwertewqn5}. Fix a
tuple $(\ell_1,\dots,\ell_n)$ of positive integers and consider
$\prod_{j=1}^n n_j^{\ell_j}$. We may assume $\ell_1 \le \cdots \le
\ell_n$ and at least one $\ell_j \ge 2$, as otherwise all
$n_j^{\ell_j}$ are prime; note there are $2^n-1$ such tuples. Using
$J_\nu(x) \ll 1$ (Lemma \ref{lem:Bessel}(\ref{lb:1})), the
contribution from this tuple is at most
\be\label{eq:addsqhigherpowers} N^{-\foh+\gep'}
\sum_{\substack{n_1^{\ell_1},\dots,n_n^{\ell_n} \le N^\sigma \\
\ell_1, \dots, \ell_r = 1; \ell_{r+1},\dots,\ell_n \ge 2}}
\frac1{\sqrt{n_1^{\ell_1}\cdots n_n^{\ell_n}}} \ \ll \
N^{-\foh+\gep'} N^{\frac{r}2\sigma}, \ee which is negligible for
$\sigma < \frac1{n-1}$ as $r \le n-1$. This completes the proof of
the lemma.
\end{proof}


\subsection{Converting from Sums to
Integrals}\label{sec:14convert}

In this subsection we prove the following lemma, which will be used
to finish the evaluation of $S_2^{(n)}$ in \S\ref{sec:14changevar}.
\begin{lem}\label{lem:doing_the_prime_sums}
Under the Riemann Hypothesis for $\zeta(s)$, if $\supp(\hphi)
\subseteq (-\frac{1}{n-1},\frac{1}{n-1})$, then
\begin{equation}\label{eq:result_doing_the_prime_sums}
\sum_{n_1,\dots,n_n} \left(\prod_{i=1}^n \hphi\left(\frac{\log
n_i}{\log R}\right) \frac{\chi_0(n_i)\Lambda(n_i)}{\sqrt{n_i}\log
R} \right) J_{k-1}\left(\frac{4\pi m \sqrt{n_1\dots
n_n}}{b\sqrt{N}}\right) \ = \ I_n(\hphi) + O\left(
N^{\frac{(n-1)\sigma}2+\gep}\right)
\end{equation}
uniformly for $m\leq N^\epsilon$ and $b\geq 1$, where
\begin{equation}\label{eq:S answer n}
I_n(\hphi)\ =\ \frac{b\sqrt{N}}{2\pi m} \int_{0}^\infty J_{k-1}(x)
\widehat{\Phi_n}\left(\frac{2\log(bx\sqrt{N} / 4\pi m)}{\log
R}\right)\frac{\d x}{\log R}
\end{equation} and $\Phi_n(x) = \phi(x)^n$. Note for $\sigma <
\frac1{n-1}$ that the main term is larger than the error term.
\end{lem}


\begin{proof}

Let $G_{k-1}(s)$ be the Mellin transform of the Bessel function. By
(6.561.14) of \cite{GR} it is
\begin{align}\label{eq:G_nu}
G_{k-1}(s) &= \int_0^\infty J_{k-1}(x) x^{s-1} \d x \nonumber\\
&= 2^{s-1} \Gamma\left(\frac{k-1+s}2\right) \Big/
\Gamma\left(\frac{k+1-s}2\right), \ \ \ \Re(s+k-1) > 0,\ \Re(s) <
\frac32.
\end{align}
Since we have $k \ge 2$, we may take $\Re(s) \in [0,1]$. The inverse
transform is
\begin{equation}\label{eq:JandGrelation}
J_{k-1}(x)\ =\  \frac1{2\pi \i} \int_{\Re(s) = 1} G_{k-1}(s) x^{-s}
\d s,
\end{equation}
and so our task is to evaluate \bea\label{eq:termMTE} & &
\sum_{n_1,\dots,n_n} \left[ \prod_{i=1}^n \hphi\left(\frac{\log
n_i}{\log R}\right) \frac{\chi_0(n_i)\Lambda(n_i)}{\sqrt{n_i} \log
R}\right] J_{k-1}\left(\frac{4\pi m \sqrt{n_1\cdots
n_n}}{b\sqrt{N}}\right) \nonumber\\ & = \ &
\sum_{n_1,\dots,n_n}\left[ \prod_{i=1}^n \hphi\left(\frac{\log
n_i}{\log R}\right) \frac{\chi_0(n_i)\Lambda(n_i)}{\sqrt{n_i} \log
R}\right] \frac1{2\pi \i} \int_{\Re(s) = 1} G_{k-1}(s)
\left(\frac{4\pi m \sqrt{n_1\cdots n_2}}{b\sqrt{N}}\right)^{-s} \d s
\nonumber\\ & = \ &
 \frac1{2\pi \i} \int_{\Re(s) = 1} \left[ \prod_{i=1}^n
\sum_{n_i}\hphi\left(\frac{\log n_i}{\log R}\right)
\frac{\chi_0(n_i)\Lambda(n_i)}{n_i^{(1+s)/2} \log R}\right]
\left(\frac{4\pi m }{b\sqrt{N}}\right)^{-s} G_{k-1}(s) \;\d s \ .
\eea Using the Mellin transform allows us to move the summations
inside the product. In \eqref{eq:tempIphizs} we derive a useful
integral version of the $n_i$-sums.

Note that for $\Re(z)>1$,
\begin{equation}\label{eq:trivial_char_L_func}
L(z,\chi_0)\ =\ \prod_{p} \left(1-\frac{\chi_0(p)}{p^z}\right)^{-1}
\ = \ \zeta(z) \prod_{p|b} \left(1-\frac{1}{p^z}\right),
\end{equation}
and so $L(z,\chi_0)$ has a simple pole at $z=1$, and zeros at $z=0$ and
all the zeros of the Riemann zeta function. Consider the integral
\begin{equation}\label{eq:defn_I}
I\ =\ -\frac{1}{2\pi\i} \int_{\Re(z)=2} \phi\left(\frac{(2z-s-1)\log
R}{4\pi\i}\right) \frac{L'}{L}(z,\chi_0)\;\d z,
\end{equation} where we have extended $\phi$ by setting \be
\phi(x+\i y) \ = \ \int_{-\infty}^\infty \hphi(u) e^{2\pi\i (x+\i y)
u}\;\d u. \ee Since $\hphi$ is a Schwartz function of compact
support, $\phi(x+\i y)$ decays rapidly as $x\to\pm\infty$ for any
fixed $y$. Thus the integral in \eqref{eq:defn_I} is absolutely
convergent, and all contour shifts are well defined. On the line of
integration the $L$-function can be written as a Dirichlet series,
and we have
\begin{align}\label{eq:tempIphizs}
I &\ = \ \frac{1}{2\pi\i} \int_{\Re(z)=2}
\phi\left(\frac{(2z-s-1)\log R}{4\pi\i}\right) \sum_{r=1}^\infty
\frac{\Lambda(r)
\chi_0(r)}{r^z}\;\d z \nonumber\\
&\ = \ \sum_{r=1}^\infty \Lambda(r) \chi_0(r) \frac{1}{2\pi\i}
\int_{\Re(z)=2} \phi\left(\frac{(2z-s-1)\log R}{4\pi\i}\right)
r^{-z}
\;\d z \nonumber\\
&\ = \ \sum_{r=1}^\infty \Lambda(r) \chi_0(r) \frac{1}{2\pi\i}
\int_{\Re(z)=(1+\Re(s))/2} \phi\left(\frac{(2z-s-1)\log
R}{4\pi\i}\right) r^{-z}
\;\d z \nonumber\\
&\ = \ \sum_{r=1}^\infty \frac{\Lambda(r) \chi_0(r)}{r^{(1+s)/2}\log
R}\ \hphi\left(\frac{\log r}{\log R}\right).
\end{align}
Interchanging of the order of summation and integration in
\eqref{eq:tempIphizs} is justified by the absolute convergence. The
$\hphi(\log r / \log R)$ factor arises from expanding the integral
in \eqref{eq:tempIphizs}; because $\Re(z) = \frac{1+\Re(s)}2$, the
argument of $\phi$ is real and the resulting integral is just the
Fourier transform.

An alternative evaluation of the integral in \eqref{eq:defn_I} is to
shift the contour to the line $\Re(z)=c$ with $1/2<c<1$. This
contour shift picks up the pole at $z=1$ and nothing else (under the
Riemann Hypothesis).

We may therefore conclude that
\begin{multline}\label{eq:prime sum answer}
\sum_{r=1}^\infty \hphi\left(\frac{\log r}{\log R}\right)
\frac{\chi_0(r)\Lambda(r)}{r^{(1+s)/2} \log R}\ =\
\phi\left(\frac{1-s}{4\pi\i}\log R\right)  \\
-\frac{1}{2\pi\i} \int_{\Re(z)=c} \phi\left(\frac{(2z-1-s)\log
R}{4\pi\i}\right) \frac{L'}{L}(z,\chi_0)\;\d z.
\end{multline}
Denoting the integral in \eqref{eq:prime sum answer} by
$\mathcal{E}(s)$, we see that \eqref{eq:termMTE} equals
\begin{multline}\label{eq:termMTE expanded}
\frac1{2\pi \i} \int_{\Re(s) = 1} \left(
\phi\left(\frac{1-s}{4\pi\i}\log R\right) + \mathcal{E}(s)
\right)^n \left(\frac{4\pi m }{b\sqrt{N}}\right)^{-s} G_{k-1}(s)
\;\d s \\
=\ \sum_{j=0}^n \binom{n}{j} \frac1{2\pi \i} \int_{\Re(s) = 1}
\phi\left(\frac{1-s}{4\pi\i}\log R\right)^{n-j} \mathcal{E}(s)^j
\left(\frac{4\pi m }{b\sqrt{N}}\right)^{-s} G_{k-1}(s) \;\d s.
\end{multline}

The main term is when $j=0$. Letting $\Phi_n(x) := \phi(x)^n$ and
using \eqref{eq:G_nu} to write $G_{k-1}(s)$ in terms of the Bessel
function, we see that it equals
\begin{align}
& \frac{1}{2\pi} \int_{-\infty}^\infty \phi\left(\frac{-t \log
R}{4\pi} \right)^n \left(\frac{4\pi m }{b\sqrt{N}}\right)^{-1-\i t}
G_{k-1}(1+\i t) \;\d t \nonumber\\ & \ = \ \frac{b \sqrt{N}}{8 \pi^2
m} \int_{-\infty}^\infty \phi\left(\frac{-t\log R}{4\pi}\right)^n
\left(\frac{4\pi m}{b\sqrt{N}}\right)^{-\i t} \int_0^\infty
J_{k-1}(x) x^{\i t} \;\d x \;\d t \nonumber\\ & \ = \
\frac{b\sqrt{N}}{8\pi^2 m } \int_0^\infty J_{k-1}(x)
\int_{-\infty}^\infty \Phi_n\left(\frac{-t\log R}{4\pi}\right)
\exp\left(\i t \log(bx\sqrt{N}/4\pi m )\right) \;\d t \;\d x
\nonumber\\ & \ = \ \frac{b\sqrt{N}}{2\pi m\log R } \int_0^\infty
J_{k-1}(x) \int_{-\infty}^\infty \Phi_n\left(u\right)
\exp\left(-2\pi\i u\ \frac{2\log(bx\sqrt{N} / 4\pi m)}{\log
R}\right) \;\d u \;\d x;
\end{align}
interchanging the order of integration is justified by Fubini's
Theorem. As the inner integral is simply the Fourier transform of
$\Phi_n$, the main term equals
\begin{equation}
\frac{b\sqrt{N}}{2\pi m \log R} \int_0^\infty J_{k-1}(x)
\widehat{\Phi_n} \left(\frac{2\log(bx\sqrt{N}/4\pi m)}{\log
R}\right) \; \d x,
\end{equation}
which we have denoted $I_n(\hphi)$.

The remaining terms in \eqref{eq:termMTE expanded} are error terms,
arising from $j \in \{1,\dots, n\}$. We shift the integrals over $s$
to the line $\Re(s)=-\gep$ and estimate $\phi$, $\mathcal{E}$ and
$G_{k-1}$ in order to bound these terms.

If $\supp \hphi \subset [-\sigma,\sigma]$ and $|x| > x_0 > 0$, then
integrating by parts $A$ times yields
\begin{align}\label{eq:bound_on_phi}
\phi\left(\frac{x+\i y}{4\pi\i} \log R\right) &\ = \ \intinf
\hphi(u)
e^{u (x+\i y) \log R / 2} \;\d u \nonumber\\
&\ = \ \intinf \hphi^{(A)}(u)
\left(\frac{-2}{(x+\i y)\log R}\right)^A e^{u (x+\i y) \log R / 2} \;\d u \nonumber\\
&\ \ll_{x_0}\ \frac{1}{(1+|y|)^A} \intinf \left|\hphi^{(A)}(u)\right| R^{x u /2} \;\d u \nonumber\\
&\ \ll_{x_0}\ \frac{N^{|x| \sigma /2}}{(1+|y|)^A}
\end{align}
since $\supp(\hphi^{(A)}) \subset [-\sigma,\sigma]$, and
$R^{|x|\sigma/2} \ll N^{|x|\sigma/2}$.

For $\Re(z)=c$ with $c \in (\foh, 1]$, we have
\begin{equation}
\frac{L'}{L}(z,\chi_0)\ \ll_c\ \log((2+|\Im(z)|) b)
\end{equation}
which follows from \eqref{eq:trivial_char_L_func} and the
well-known bound $\frac{\zeta'}{\zeta}(z) \ll_c \log(2+|\Im(z)|)$
(see, for example, Theorem 14.5 of \cite{Tit}). Therefore, if
$s=-\gep+\i t$ and $c=1/2+\gep'$, we have
\begin{align}\label{eq:bound_on_E}
\mathcal{E}(-\gep+\i t) &\ = \  -\frac{1}{2\pi\i}
\int_{\Re(z)=1/2+\gep'} \phi\left(\frac{(2z-1+\gep-\i t)\log
R}{4\pi\i}\right) \frac{L'}{L}(z,\chi_0)\;\d
z \nonumber\\
&\ \ll\ \int_{-\infty}^\infty \left|\phi\left(\frac{(2\gep'+2\i
y+\gep-\i t)\log R}{4\pi\i}\right) \right|\ \log((2+|y|) b) \;\d y \nonumber\\
&\ \ll\  \int_{-\infty}^\infty  \frac{N^{(2\gep'+\gep) \sigma
/2}}{(1+|2y - t|)^A}\ \log((2+|y|) b) \;\d y \nonumber\\
&\ \ll\ N^{\gep''} \log((2+|t|)b);
\end{align} above we used $2\gep'+\gep > 0$ so as to be able to
apply the bounds from \eqref{eq:bound_on_phi}.

We also need an estimate for $G_{k-1}$. From \eqref{eq:G_nu} and
\begin{equation}\label{eq:gammalargetexp} \Gamma(\sigma + \i t) \ =
\ \sqrt{2\pi}\ (\i t)^{\sigma - \foh}\ e^{-\frac{\pi}2|t|}\
|t/e|^{\i t}\ \left(1 + O\left(|t|^{-1}\right)\right),
\end{equation}
we have
\begin{equation}\label{eq:bound_on_G}
\left|G_{k-1}(-\gep+\i t)\right|\ =\ 2^{-\gep-1} \left|
\frac{\Gamma((k-1-\gep+\i t)/2)}{\Gamma((k-1-\gep-\i t)/2 + 1+\gep)}
\right|\ \ll\ \frac{1}{(1+|t|)^{1+\gep}} .
\end{equation}

Using \eqref{eq:bound_on_phi}, \eqref{eq:bound_on_E} and
\eqref{eq:bound_on_G} in \eqref{eq:termMTE expanded}, we have
\begin{align}
\frac1{2\pi \i} &\int_{\Re(s) = -\gep}
\phi\left(\frac{1-s}{4\pi\i}\log R\right)^{n-j} \mathcal{E}(s)^j
\left(\frac{4\pi m
}{b\sqrt{N}}\right)^{-s} G_{k-1}(s) \;\d s \nonumber\\
&\ll\ \left(\frac{m}{b\sqrt{N}}\right)^{\gep} \intinf
\left|\phi\left(\frac{1+\gep-\i t}{4\pi\i}\log R\right)\right|^{n-j}
\left|\mathcal{E}(-\gep+\i t)\right|^j \left| G_{k-1}(-\gep+\i
t)\right|
\;\d t \nonumber\\
&\ll\  \left(\frac{m}{b\sqrt{N}}\right)^{\gep} \intinf
\left(\frac{N^{(1+\gep) \sigma /2}}{(1+|t|)^A}\right)^{n-j}
\left(N^{\gep''}\log((2+|t|)b)\right)^j \frac{1}{(1+|t|)^{1+\gep}}  \;\d t \nonumber\\
&\ll\  N^{(n-j) \sigma /2 + \gep'''}.
\end{align}
Note the $t$-integral converges (it is only when $j=n$ that we need
to use $\gep>0$ to ensure convergence). The worst term is clearly
when $j=1$, and this yields the desired error term.

This completes the proof of Lemma~\ref{lem:doing_the_prime_sums}.
\end{proof}

\begin{rek} In \S\ref{sec:14changevar} we finish the evaluation of
$S_2^{(n)}$. We multiply our terms by $N^{-1/2}$ and execute the
summations over $b$ and $m$. Thus in order for the error term in
Lemma \ref{lem:doing_the_prime_sums} to be negligible we need \be
N^{-\foh} N^{\frac{(n-1)\sigma}2+\gep} \ \ll \ N^{-\gep'}, \ee which
forces $\sigma < \frac1{n-1}$. We thus see that, in the number
theory calculations, $\frac1{n-1}$ is a real boundary for this
method when we split by sign. This is very different than related
problems in \cite{Ru,Gao} and the non-split case of Theorem
\ref{thm:nosplitmoments}; the reason is due to support problems when
$n \ge 2$ from handling the Bessel-Kloosterman terms from
$\gl_f(N)$. Thus when we split by sign, we expect our results for
support up to $\frac1{n-1}$ and not $\frac2n$.
\end{rek}


\subsection{Evaluating $S_2^{(n)}$}\label{sec:14changevar}

We finish the proof of Theorem \ref{thm:extended moments} by
completing the evaluation of $S_2^{(n)}$.

\begin{lem}\label{lem:final answer for S_2}
Under GRH for $L(s,f)$ and for all Dirichlet $L$-functions, if
$n\geq 2$ and $\supp(\hphi) \subset (-\frac1{n-1},\frac1{n-1})$,
then
\begin{equation}
S_2^{(n)}\ =\ 2^{n-1} \left[ \int_{-\infty}^\infty \phi(x)^n
\frac{\sin 2\pi x}{2\pi x} \;\d x- \frac12 \phi(0)^n \right] +
O\left( \frac{k \log\log kN}{\log kN}\right).
\end{equation}
\end{lem}

\begin{proof}
Combining Lemmas \ref{lem:expand_S_2_with_prim_char} and
\ref{lem:doing_the_prime_sums}, we have shown that under GRH for
$L(s,f)$ and for all Dirichlet $L$-functions, if $n\geq 2$ and
$\supp(\hphi) \subset (-\frac1{n-1},\frac1{n-1})$, then
\begin{equation}\label{eq:lemncmtodate}
S_2^{(n)}\ =\  - \frac{2^{n+1}\pi}{\sqrt{N}} \sum_{m \le N^\epsilon}
\frac{1}{m}\sum_{\substack{(b,N) \\ b < N^{2006}}}
\frac{R(m^2,b)R(1,b)}{b\varphi(b)} \left(I_n(\hphi) + O\left(
N^{\frac{(n-1)\sigma}2+\gep}\right) \right)
\end{equation}
where
\begin{equation}\label{eq:sn2inhphi}
I_n(\hphi) \ = \ \frac{b\sqrt{N}}{2\pi m} \int_{x=0}^\infty
J_{k-1}(x) \widehat{\Phi_n}\left(\frac{2\log(bx\sqrt{N} / 4\pi
m)}{\log R}\right)\frac{\d x}{\log R}
\end{equation} and $\Phi_n(x) = \phi(x)^n$.
Since by \eqref{eq:defn Rnew} we have $R(m^2,b)R(1,b) \ll m^4$, the
contribution from the $O$-term in \eqref{eq:lemncmtodate} is bounded
by \bea N^{-\foh} \sum_{m \le N^\gep} \frac{m^4}{m}
\sum_{\substack{(b,N) \\ b < N^{2006}}} \frac1{b^2}\left(
N^{\frac{(n-1)\sigma}2+\gep}\right)\ \ll \ N^{\frac{n-1}2 \cdot
(\sigma - \frac1{n-1} + \gep'')}, \eea which is $O(N^{-\gep'''})$
for $\sigma < \frac1{n-1}$.

We are left with evaluating the main term. The rapid decay of
$I_n(\hphi)$ with respect to $b$ allows us to extend the $b$-sum of
the main term of \eqref{eq:lemncmtodate} to all $b$ relatively prime
to $N$. From \eqref{eq:sn2inhphi} we have \be I_n(\hphi) \ \ll \
\frac{b\sqrt{N}}{2\pi m}  \int_0^\infty
\left|\widehat{\Phi_n}\left(\frac{2\log(u\sqrt{N} / 4\pi m)}{\log
R}\right)\right| \frac{\d u}{b} \ \ll \ \frac{\sqrt{N}}{m}. \ee From
\eqref{eq:defn Rnew} we have $R(m^2,b) R(1,b) \ll m^4$. The $m$-sum
is $O(N^{4\gep})$, the factor of $N^{-\foh}$ cancels the factor of
$N^{\foh}$, and we have a $b$-sum of $b^{-2}$ (which is negligible
for the terms with $b \ge N^{2006}$).

As $\Phi_j(x) = \phi(x)^j$, $\widehat{\Phi_n}$ is the convolution of
$\phi$ with itself $n$ times. In particular we have
\begin{equation}
\widehat{\Phi_n}(u) \ =\ \int_{-\infty}^\infty
\widehat{\Phi_{n-1}}(w) \hphi(u-w) \;\d w.
\end{equation}
Note that the support of $\widehat{\Phi_n}$ is at most $n$ times
that of $\hphi$, which means for $n\ge 2$ it is less than
$\frac{n}{n - 1} \le 2$. Therefore we may apply
Lemma~\ref{lem:ilschap7}. We find that the main term of $S_2^{(n)}$
is
\begin{align}\label{eq:5.19}
& - \frac{2^{n+1}\pi}{\sqrt{N}}\sum_{m \le N^\epsilon}
\frac{1}{m}\sum_{(b,N) = 1} \frac{R(m^2,b)R(1,b)}{b\varphi(b)}\
I_n(\hphi)\nonumber\\  &
\quad = \ -
\frac{2^{n+1}\pi}{\sqrt{N}}\sum_{m \le N^\epsilon} \frac{1}{m}
\sum_{(b,N) = 1} \frac{R(m^2,b) R(1,b)}{b\varphi(b)}
\frac{b\sqrt{N}}{2\pi m} \int_{0}^\infty J_{k-1}(y)
\widehat{\Phi_n}\left( \frac{2 \log (by\sqrt{N} / 4\pi m) }{\log
R}\right) \frac{\d y}{\log R}\nonumber\\
& \quad = \ -
\frac{2^{n+1}\pi}{\sqrt{N}}\sum_{m \le N^\epsilon} \frac{1}{m}
\sum_{(b,N) = 1} \frac{R(m^2,b) R(1,b)}{b\varphi(b)}
\frac{b\sqrt{N}}{2\pi m} \int_{0}^\infty J_{k-1}(y)
\widehat{\Phi_n}\left( \frac{2 \log (by\sqrt{N} / 4\pi m) }{\log
R}\right) \frac{\d y}{\log R}\nonumber\\
& \quad = \ -2^n \cdot
\left(-\frac12\right) \cdot \left[ \int_{-\infty}^\infty \Phi_{n}(x)
\frac{\sin 2\pi x}{2\pi x} \;\d x- \frac12 \Phi_{n}(0) \right]
+O\left(
\frac{k \log\log kN}{\log kN}\right)\nonumber\\
& \quad = \ 2^{n-1} \left[ \int_{-\infty}^\infty \phi(x)^n \frac{\sin
2\pi x}{2\pi x} \;\d x- \frac12 \phi(0)^n \right]+O\left( \frac{k
\log\log kN}{\log kN}\right).
\end{align}
This completes the proof of the lemma, and Theorem \ref{thm:extended
moments}.
\end{proof}

\begin{rek}
Note that if $\supp(\hphi)  \subset (-\frac1{n}, \frac{1}{n})$ then
$\supp(\widehat{\Phi_n})  \subset (-1,1)$. In this region the
Kloosterman-Bessel terms are negligible, and the contribution to the
centered moment in \eqref{eq:nth_mmt_in_terms_of_S} from $S_2^{(n)}$
vanishes. As $\frac{1}{n-1}> \frac1{n}$ for $n\ge 2$, we have
entered the non-trivial region where these contributions do not
vanish. Thus the mock Gaussian result of
Theorem~\ref{thm:mock-Gaussian for D} is sharp.
\end{rek}



\section{Random Matrix Theory: Proof of Theorem~\ref{thm:mmts_RMT}}
\label{sect:RMT moments}

In this section we prove Theorem~\ref{thm:mmts_RMT} by calculating
the centered moments of $\Zf(U)$ when averaged over $\soe$ and
$\soo$. For small support the moments agree with those of the
Gaussian; for larger support, however, the moments differ.
\subsection{Introduction}

Let $U$ be an $M\times M$ unitary matrix with eigenvalues
$e^{\i\gt_1}$, \dots, $e^{\i\gt_M}$. For a real, even integrable
function $\phi$ which decays sufficiently rapidly, define
\begin{align}
F_M(\gt) &\ = \ \sum_{j=-\infty}^\infty
\phi\left(\frac{M}{2\pi}(\theta + 2\pi j) \right)\nonumber\\
&\ =\ \frac{1}{M} \sum_{k=-\infty}^\infty
\hphi\left(\frac{k}{M}\right) e^{\i k \theta},
\end{align}
which is a $2\pi$-periodic function emphasizing points near $0
\mod 2\pi$. Define
\begin{equation}
\Zf(U)\ =\ \sum_{n=1}^M F_M\left(\theta_n\right).
\end{equation}
This is the random matrix equivalent of $D(f;\phi)$. More
precisely, moments of $D(f,\phi)$ averaged over $f\in \hkpn$ as
$N\to\infty$ should correspond to moments of $\Zf(U)$ when
averaged with respect to Haar measure over $\SO(M)$ matrices as
$M$ tends to infinity through \emph{even} integers, while moments
of $D(f,\phi)$ averaged over $f\in \hkmn$ as $N\to\infty$ should
correspond to moments of $\Zf(U)$ when averaged with respect to
Haar measure over $\SO(M)$ matrices as $M$ tends to infinity
through \emph{odd} integers.
\begin{rek}
If we restrict the eigenangles such that $-\pi < \theta_n \leq
\pi$, then
\begin{equation}
\Zf(U)\ \sim\ \sum_{n=1}^M
\phi\left(\frac{M}{2\pi}\theta_n\right).
\end{equation}
However, using $F_M(\theta)$ in the definition of $\Zf(U)$ is more
natural because the eigenangles of orthogonal matrices are
$2\pi$-periodic.
\end{rek}

Much of the work required to calculate the moments of $\Zf(U)$ was
done in the paper of Hughes and Rudnick \cite{HR1} (building on work
of Soshnikov \cite{Sosh}), and we simply quote the results we need
to show Theorem~\ref{thm:mmts_RMT}. The novelty here is
desymmetrizing the integrals to handle the combinatorics in the
non-trivial range. This is necessary in order to write the formulas
in such a way as to facilitate comparisons with number theory.

Theorem~5 of \cite{HR1}, when applied to the case $\Zf(U)$, shows
that the means over $\soe$ and $\soo$ are
\begin{align}
C_1^{\soe} \ &:=\ \lim_{\substack{M\to\infty \\ M {\rm even}}}
\E_{\SO(M)} \left[\Zf(U)
\right] \ =\ \hphi(0) + \frac12\int_{-1}^1 \hphi(y)\;\d y \\
C_1^{\soo} \ &:=\ \lim_{\substack{M\to\infty \\ M {\rm odd}}}
\E_{\SO(M)} \left[\Zf(U) \right] \ =\ \hphi(0) + \frac12\int_{-1}^1
\hphi(u)\;\d u + \int_{|y|\geq 1} \hphi(u)\;\d u,
\end{align}
respectively, where $\E_{\SO(M)}$ denotes expectation with respect
to Haar measure over the classical compact group of $M\times M$
special orthogonal matrices. Furthermore, that theorem states that
the variance of $\Zf(U)$ over $\soe$ is \bea\label{eq:variance
SO(even)} C_2^{\soe} &  := & \lim_{\substack{M\to\infty \\ M {\rm
even}}} \E_{\SO(M)}
\left[\left(\Zf(U) - C_1^{\soe}\right)^2 \right]  \nonumber\\
&=& 2\int_{-\infty}^\infty \min(|y|,1) \hphi(y)^2 \;\d y +
2\int_{y=-1/2}^{1/2} \int_{|x|\geq 1/2} \hphi(x+y) \hphi(x-y) \;\d
x\;\d y,\nonumber\\ \eea and over $\soo$ is \bea\label{eq:variance
SO(odd)} C_2^{\soo} \ &:= & \ \lim_{\substack{M\to\infty \\ M {\rm
odd}}} \E_{\SO(M)}
\left[\left(\Zf(U) - C_1^{\soo}\right)^2 \right] \nonumber\\
&=& 2\int_{-\infty}^\infty \min(|y|,1) \hphi(y)^2 \;\d y -
2\int_{y=-1/2}^{1/2} \int_{|x|\geq 1/2} \hphi(x+y) \hphi(x-y) \;\d
x\;\d y.\nonumber\\ \eea

Changing variables to $u=x+y$ and $v=x-y$ we see that
\begin{equation}
\int_{y=-1/2}^{1/2} \int_{|x|\geq 1/2} \hphi(x+y) \hphi(x-y) \;\d
x\;\d y\ =\  \frac12\iint_{A} \hphi(u) \hphi(v)  \;\d u \;\d v,
\end{equation}
where
\begin{equation}
A\ =\ \{ |u+v| \geq 1 \}\cap\{|u-v|\leq 1\} .
\end{equation}
Note that if $|u| \leq 1$ and $|v| \leq 1$, then whenever $|u+v|
\geq 1$ we have $\{|u-v|\leq 1\}$. Therefore if $\supp(\hphi)
\subseteq [-1,1]$, we have \bea\label{eq:neededforn=2}
\frac12\iint_{A} \hphi(u) \hphi(v)  \;\d u \;\d v & = &
\frac12\int_{-1}^1 \int_{-1}^1 \hphi(u) \hphi(v)
\I_{\{|u+v|\geq 1\}} \;\d u \;\d v\nonumber\\
&=&-\left(\int_{-\infty}^\infty \phi(x)^2 \frac{\sin 2\pi x}{2 \pi
x} \;\d x - \frac12 \phi(0)^2\right).\ \ \ \
\label{eq:RMT_variance_restricted} \eea the last line following from
the Fourier transform identity (see Lemma \ref{thm:neededn=2int} for
a proof)
\begin{align}
& \int_{-\infty}^\infty \phi(x)^2 \frac{\sin 2\pi x}{2 \pi x} \;\d
x  - \frac12 \phi(0)^2\nonumber\\ &\ \ \ = \  \frac12 \iint
\hphi(u) \hphi(v) \I_{\{|u+v|\leq 1\}} \;\d u \;\d v - \frac12
\iint \hphi(u) \hphi(v) \;\d u\;\d v \nonumber\\
&\ \ \  = \  -\frac12 \iint \hphi(u) \hphi(v) \I_{\{|u+v|
> 1\}} \;\d u \;\d v .
\end{align}
Furthermore, note that if either $|u|>1$ or $|v|>1$ then $|u+v|
\geq 1$ does not necessarily imply $|u-v|\leq 1$, and so
\eqref{eq:RMT_variance_restricted} does not hold if the support of
$\hphi$ exceeds $[-1,1]$.

This proves Theorem~\ref{thm:mmts_RMT} in the case $n=1$ and
$n=2$. While the higher moments of $\Zf(U)$ can be calculated
using Weyl's explicit representation of Haar measure for even and
odd orthogonal groups, as in \cite{HR1} we deal with its
cumulants. Denote
\begin{equation}
\sum_{\ell=1}^\infty C_\ell^{\soe} \frac{\lambda^\ell}{\ell!}\ =\
\lim_{\substack{M\to\infty\\M {\rm even}}} \log \E_{\SO(M)} \left[
\exp\left(\lambda \Zf(U)\right) \right]
\end{equation}
and
\begin{equation}
\sum_{\ell=1}^\infty  C_\ell^{\soo} \frac{\lambda^\ell}{\ell!} \ =
\ \lim_{\substack{M\to\infty\\M {\rm odd}}} \log \E_{\SO(M)}
\left[ \exp\left(\lambda \Zf(U)\right) \right].
\end{equation}
Knowing the first $n$ cumulants is equivalent to knowing the first
$n$ moments, which is evident from the identity
\begin{equation}\label{eq:mmts in terms of cumulants}
\E_{\SO(M)}\left[\left(\Zf(U)\right)^n \right]\ =\ \sum
\left(\frac{C_1}{1!}\right)^{k_1}\left(\frac{C_2}{2!}\right)^{k_2}\cdots\left(\frac{C_n}{n!}\right)^{k_n}
\frac{n!}{k_1!\;k_2!\cdots k_n!},
\end{equation}
where the sum runs over all non-negative values of $k_j$
($j=1,\dots,n$) such that $\sum_{j=1}^n j k_j = n$.
Theorem~\ref{thm:mock_Gaussian} implies that if $j\geq 3$ and
$\supp(\hphi) \subseteq [-\frac{1}{j},\frac{1}{j}]$, then both
$C_j^{\soe} = 0$ and $C_j^{\soo} = 0$. Therefore, by
\eqref{eq:mmts in terms of cumulants}, if $\supp(\hphi) \subseteq
[-\frac{1}{n-1} , \frac{1}{n-1}]$ for $n\geq 3$, the only terms
which contribute to the $n$\textsuperscript{th} moment are $C_1$,
$C_2$ and $C_n$, and we have
\begin{equation}\label{eq:cmcumC}
\lim_{\substack{M\to\infty\\M {\rm even}}} \left[ \left(\Zf(U) -
C_1^{\soe} \right)^n \right] =
\begin{cases}
C_{2k}^{\soe} + \left(C_2^{\soe}\right)^k \dfrac{(2k)!}{2^k k!} &
\text{$n=2k$}\\
C_{2k+1}^{\soe} & \text{$n=2k+1$,}
\end{cases}
\end{equation}
and similarly for $\soo$. If $\supp(\hphi) \subseteq [-\foh, \foh]$,
\eqref{eq:variance SO(even)} and \eqref{eq:variance SO(odd)} yield
\begin{equation}
C_2^{\soe}\ =\ C_2^{\soo}\ =\ 2\int_{-1/2}^{1/2} |y| \hphi(y)^2\;\d
y \ = \ \sigma_\phi^2,
\end{equation}
where $\sigma_\phi^2$ is given by \eqref{eq:defRnn-1b}. Therefore,
by \eqref{eq:cmcumC}, Theorem~\ref{thm:mmts_RMT} will follow from
showing
\begin{equation}\label{eq:RTP cumulant even}
C_n^{\soe}\ =\ (-1)^{n-1} 2^{n-1}\left[ \int_{-\infty}^\infty
\phi(x)^{n} \frac{\sin 2\pi x}{2\pi x} \;\d x - \foh \phi(0)^{n}
\right]
\end{equation}
and
\begin{equation}\label{eq:RTP cumulant odd}
C_n^{\soo}\ =\ (-1)^n 2^{n-1}\left[ \int_{-\infty}^\infty
\phi(x)^{n} \frac{\sin 2\pi x}{2\pi x} \;\d x - \foh \phi(0)^{n}
\right]
\end{equation}
for $n\geq 3$ and $\supp(\hphi) \subseteq [-\frac1{n-1} ,
\frac1{n-1}]$.

Let
\begin{multline}\label{eq:defn R_n}
Q_n(\phi) := 2^{n-1} \sum_{m=1}^n
\sum_{\substack{\lambda_1+\cdots+\lambda_m = n\\ \lambda_j \geq
1}} \frac{(-1)^{m+1}}{m} \frac{n!}{\lambda_1! \cdots \lambda_m!}
\intinf \cdots \intinf
\phi(x_1)^{\lambda_1} \cdots \phi(x_m)^{\lambda_m} \\
\times S(x_1-x_2) S(x_2-x_3) \cdots S(x_{m-1}-x_m) S(x_m + x_1) \;\d
x_1 \cdots \d x_m,
\end{multline}
where
\begin{equation}
S(x)\ =\ \frac{\sin(\pi x)}{\pi x}\ =\ \intinf \I_{\{|u|\leq 1/2\}}
e^{2\pi\i x u} \;\d u.
\end{equation}
To prove (\ref{eq:RTP cumulant even}, \ref{eq:RTP cumulant odd}), we
again use results from \cite{HR1} (Section 2.1, Lemma 6 and Theorem
7), where it was shown that if $\supp(\hphi) \subseteq
[-\frac2{n},\frac2{n}]$, then $C_{n}^{\soe} = Q_n(\phi)$ and
$C_{n}^{\soo} = -Q_n(\phi)$. Therefore we must show that whenever
$\supp(\hphi) \subseteq [-\frac{1}{n-1} , \frac{1}{n-1}]$,
$Q_n(\phi)$, defined in \eqref{eq:defn R_n}, can also be written as
\bea\label{eq:eqforRn} Q_n(\phi) & \ = \ & (-1)^n 2^{n-2}
\idotsint_{-\infty}^\infty \hphi(u_1)\cdots \hphi(u_n)
\I_{\{|u_1+\dots +u_n|\geq 1\}} \;\d u_1 \cdots \d u_n\nonumber\\
&= &(-1)^{n-1} 2^{n-1}\left[ \int_{-\infty}^\infty \phi(x)^{n}
\frac{\sin 2\pi x}{2\pi x} \;\d x - \foh \phi(0)^{n} \right]; \eea
the two expressions for $Q_n(\phi)$ in \eqref{eq:eqforRn} are
equal by Lemma \ref{eq:FT id1}.

To prove \eqref{eq:eqforRn}, we use Plancherel's identity in
\eqref{eq:defn R_n}, and write the test function $\phi$ in terms
of its Fourier transform $\hphi$, and $S(x)$ in terms of its
Fourier transform, obtaining
\begin{multline}
Q_n(\phi)\ =\ 2^{n-1} \sum_{m=1}^n
\sum_{\substack{\lambda_1+\cdots+\lambda_m = n\\ \lambda_j \geq
1}} \frac{(-1)^{m+1}}{m} \frac{n!}{\lambda_1! \cdots
\lambda_m!} \intinf \cdots \intinf \hphi(y_1) \cdots \hphi(y_n) \\
e^{2\pi\i x_1 \left(u_1 + u_m + y_1 + \cdots +y_{\lambda_1}
\right) } e^{2\pi\i x_2 \left(u_2 - u_1
+y_{\lambda_1+1}+\cdots+y_{\lambda_1+\lambda_2} \right) } \cdots
e^{2\pi\i x_m \left(u_m - u_{m-1} + y_{\lambda_1+\cdots+\lambda_{m-1}+1}+\cdots+y_n \right) }\\
\times \I_{\{|u_1|\leq 1/2\}} \cdots \I_{\{|u_m|\leq 1/2\}} \;\d u_1
\cdots \d u_m \;\d y_1 \cdots \d y_n \;\d x_1 \cdots \d x_m.
\end{multline}
For simplicity write \bea
Y_1 &\ := \ & y_1 + \cdots + y_{\lambda_1}\nonumber\\
Y_2 &:= & y_{\lambda_1+1}+\cdots+y_{\lambda_1+\lambda_2}\nonumber\\
& \vdots & \nonumber\\ Y_m &:= &
y_{\lambda_1+\cdots+\lambda_{m-1}+1}+\cdots+y_n. \eea Integrating
over $x_1$ to $x_m$ converts the exponentials to delta
functionals, and we get \bea Q_n(\phi) & \ = \ & 2^{n-1}
\sum_{m=1}^n \sum_{\substack{\lambda_1+\cdots+\lambda_m = n\\
\lambda_j \geq 1}} \frac{(-1)^{m+1}}{m} \frac{n!}{\lambda_1!
\cdots
\lambda_m!} \intinf \cdots \intinf \hphi(y_1) \cdots \hphi(y_n) \nonumber\\
& &  \ \ \times\ \delta \left(u_1 + u_m + Y_1 \right) \delta
\left(u_2 - u_1 +Y_2 \right)  \cdots
\delta \left(u_m - u_{m-1} + Y_m \right) \nonumber\\
& & \ \ \times\ \I_{\{|u_1|\leq 1/2\}} \cdots \I_{\{|u_m|\leq
1/2\}} \;\d u_1 \cdots \d u_m \;\d y_1 \cdots \d y_n. \eea
Changing variables to
\begin{align}
v_1 &\ :=\ u_1+u_m & u_1 &\ =\ \frac12(v_1-v_2-\dots-v_m)\nonumber\\
v_2 &\ :=\ u_2-u_1 & u_2 &\ =\ \frac12(v_1+v_2-v_3-\dots-v_m)\nonumber\\
&\ \ \ \vdots\ & &\ \ \vdots\ \nonumber\\
v_m &\ :=\ u_m - u_{m-1} & v_m &\ = \ \frac12(v_1+v_2+\cdots+v_m)
\end{align}
(the Jacobian from this transformation is $\foh$) leads to
\begin{multline}
Q_n(\phi) = 2^{n-1} \sum_{m=1}^n
\sum_{\substack{\lambda_1+\cdots+\lambda_m = n\\ \lambda_j \geq
1}} \frac{(-1)^{m+1}}{m} \frac{n!}{\lambda_1! \cdots
\lambda_m!} \intinf \cdots \intinf \hphi(y_1) \cdots \hphi(y_n) \\
\times \frac12 \I_{\{|Y_1-Y_2-Y_3-\cdots-Y_m|\leq 1\}}
\I_{\{|Y_1+Y_2-Y_3-\cdots-Y_m|\leq 1\}} \cdots
\I_{\{|Y_1+Y_2+Y_3+\cdots+Y_m|\leq 1\}}  \;\d y_1 \cdots \d y_n.
\end{multline}
Making use of the fact that $\hphi$ is an even function, we
desymmetrize this by writing
\begin{equation}
Q_n(\phi)\ =\ 2^{n-2} \int_0^\infty \cdots \int_0^\infty
\hphi(y_1) \cdots \hphi(y_n) K(y_1,\dots,y_n) \;\d y_1 \cdots \d
y_n ,
\end{equation}
where
\begin{equation}
K(y_1,\dots,y_n) =  \sum_{m=1}^n \sum_{\substack{\lambda_1+\cdots+\lambda_m = n\\
\lambda_j \geq 1}} \frac{(-1)^{m+1}}{m} \frac{n!}{\lambda_1!
\cdots \lambda_m!} \sum_{\epsilon_1=\pm 1 , \dots , \epsilon_n =
\pm 1} \prod_{\ell=1}^m \I_{\{\left|\sum_{j=1}^n \eta(\ell,j)
\epsilon_j y_j\right| \leq 1 \}}
\end{equation}
with
\begin{equation}
\eta(\ell,j)\ =\
\begin{cases}
+1 & \text{ if } j \leq \sum_{k=1}^\ell \lambda_k\\
-1 & \text{ if } j >  \sum_{k=1}^\ell \lambda_k.
\end{cases}
\end{equation}

If $0\leq y_j \leq \frac{1}{n}$ for all $j$, then
\begin{equation}
\prod_{\ell=1}^m \I_{\{\left|\sum_{j=1}^n \eta(\ell,j) \epsilon_j
y_j\right| \leq 1 \}}\ =\ 1
\end{equation}
for all choices of $\epsilon_j=\pm 1$. There are $2^n$ choices of
possible $n$-tuples $(\gep_1,\dots,\gep_n)$, and so if $n\geq 2$,
\begin{equation}\label{eq:K small y}
K(y_1,\dots,y_n)\ =\  \sum_{m=1}^n \sum_{\substack{\lambda_1+\cdots+\lambda_m = n\\
\lambda_j \geq 1}} \frac{(-1)^{m+1}}{m} \frac{n!}{\lambda_1!
\cdots \lambda_m!} 2^n\ =\ 0,
\end{equation}
which follows from a trick of Soshnikov \cite{Sosh} obtained by
evaluating the generating series
\begin{equation}\label{eq:gen_func_1}
z\ =\ \log(1+(e^z-1))\ =\ \sum_{n=1}^\infty z^n\ \sum_{m=1}^n\
\sum_{\substack{\lambda_1+\cdots+\lambda_m = n\\
\lambda_j \geq 1}}\ \frac{(-1)^{m+1}}{m} \frac{1}{\lambda_1! \cdots
\lambda_m!}.
\end{equation}
If $0\leq y_j \leq \frac{1}{n-1}$ for all $j$, and $y_1 + \cdots +
y_n \geq 1$ (so at least one $y_j > \frac{1}{n}$), then
\begin{equation}
\I_{\{\left|\sum_{j=1}^n \eta(\ell,j) \epsilon_j y_j\right| \leq 1
\}}\ =\ 0
\end{equation}
if and only if either $\eta(\ell,j) \epsilon_j = +1$ for all $j$,
or $\eta(\ell,j) \epsilon_j = -1$ for all $j$. Note there are
exactly $2m$ choices for the $n$-tuple $(\epsilon_1, \dots,
\epsilon_n)$ which  yield
\begin{equation}\label{eq:eptuple}
\prod_{\ell=1}^m \I_{\{\left|\sum_{j=1}^n \eta(\ell,j) \epsilon_j
y_j\right| \leq 1 \}}\ =\ 0,
\end{equation}
and the remaining $2^n- 2m$ choices yield the product equals $1$.
This follows from the $\eta(\ell,j)$ change signs so that no
choice of $(\epsilon_1, \dots, \epsilon_n)$ makes two terms in the
product vanish. There are $m$ factors in the product, and each
factor is zero for exactly two choices of $(\epsilon_1, \dots,
\epsilon_n)$.

Hence for $n\geq 2$,
\begin{equation}\label{eq:K big y}
K(y_1,\dots,y_n)\ =\  \sum_{m=1}^n \sum_{\substack{\lambda_1+\cdots+\lambda_m = n\\
\lambda_j \geq 1}} \frac{(-1)^{m+1}}{m} \frac{n!}{\lambda_1!
\cdots \lambda_m!} \left(2^n-2m\right)\ =\ 2(-1)^n,
\end{equation}
which comes from evaluating the coefficient of $z^n$ in
\eqref{eq:gen_func_1} and in the generating series
\begin{equation}
\sum_{n=0}^\infty \frac{(-1)^n}{n!} z^n = e^{-z}\ =
\ \frac{1}{1+(e^z-1)}\ =\ \sum_{n=1}^\infty z^n \sum_{m=1}^n \sum_{\substack{\lambda_1+\cdots+\lambda_m = n\\
\lambda_j \geq 1}} (-1)^{m} \frac{1}{\lambda_1! \cdots
\lambda_m!}.
\end{equation}
Combining \eqref{eq:K small y} and \eqref{eq:K big y}, we see that
if $\supp(\hphi) \subseteq [-\frac{1}{n-1} , \frac{1}{n-1}]$, then
\begin{equation}
Q_n(\phi)\ =\ (-1)^n 2^{n-1} \int_0^{\frac{1}{n-1}} \cdots
\int_0^{\frac{1}{n-1}} \hphi(y_1) \cdots \hphi(y_n) \I_{\{y_1 +
\cdots + y_n \geq 1\}} \;\d y_1 \cdots \d y_n.
\end{equation}
The final step in the proof of Theorem \ref{thm:mmts_RMT} is the
observation (see Lemma \ref{thm:neededn=2int}) that
\begin{multline}
\intinf \phi(x)^n \frac{\sin 2\pi x}{2\pi x} \;\d x -
\frac12\phi(0)^n \\
=\ \frac12 \intinf\cdots\intinf \hphi(y_1) \cdots \hphi(y_n)
\left(\I_{\{|y_1 + \cdots + y_n| \leq 1\}} - 1\right) \;\d y_1
\cdots \;\d y_n.
\end{multline}
Furthermore, if we assume that $\supp(\hphi) \subseteq
[-\frac{1}{n-1} , \frac{1}{n-1}]$, then $\left(\I_{\{|y_1 + \cdots
+ y_n| \leq 1\}} - 1\right)$ equals zero if the $y_j$ are not all
of the same sign. Under this assumption, we therefore may write
\begin{align}
& \frac12 \intinf\cdots\intinf \hphi(y_1) \cdots \hphi(y_n)
\left(\I_{\{|y_1 + \cdots + y_n| \leq 1\}} - 1\right) \;\d y_1 \cdots \;\d y_n \nonumber\\
&\ \ \ \ =\ \int_0^{\frac{1}{n-1}} \cdots \int_0^{\frac1{n-1}}
\hphi(y_1) \cdots \hphi(y_n) \left(\I_{\{y_1 + \cdots + y_n \leq
1\}} - 1\right) \;\d
y_1 \cdots \;\d y_n \nonumber\\
&\ \ \ \ =\ -\int_0^{\frac{1}{n-1}} \cdots \int_0^{\frac{1}{n-1}}
\hphi(y_1) \cdots \hphi(y_n) \I_{\{y_1 + \cdots + y_n > 1\}} \;\d
y_1 \cdots \d y_n.
\end{align}

We therefore conclude that if $\supp(\hphi) \subseteq
[-\frac{1}{n-1}, \frac{1}{n-1}]$, then
\begin{equation}\label{eq:6.41}
Q_n(\phi)\ =\ (-1)^{n-1} 2^{n-1} \left(\intinf \phi(x)^n
\frac{\sin 2\pi x}{2\pi x} \;\d x - \frac12\phi(0)^n \right)
\end{equation}
as required.

\begin{rek} Note that $[-\frac1{n-1},\frac1{n-1}]$ is a natural
boundary. We crucially used each $y_i \le \frac1{n-1}$ in showing
there are exactly $2m$ choices for the $\gep$-tuples which make
\eqref{eq:eptuple} vanish. Indeed, beyond this point the kernel does
not have the shape of \eqref{eq:6.41}, indicating the presence of
additional terms. On the number theory side, these terms will arise
from a more detailed study of the prime powers in
Lemma~\ref{lem:natural boundary}. The new terms cannot arise from
the integral in \eqref{eq:5.19}, as this hold for $\hphi$ supported
up to $(-\frac2n,\frac2n)$.
\end{rek}


\section{The order of vanishing of $L$-functions at the critical
point}\label{sec:ordervanLFn}

We show how Corollary~\ref{cor:boundsordvancp} follows from
Theorem~\ref{thm:extended moments}. We need to assume GRH for
$L(s,f)$, which means that all non-trivial zeros are on the critical
line; this allows us to obtain bounds for the number of zeros at the
central point by using non-negative test functions. Note this rate
of decay is significantly better than previous estimates.

We compare our results to the bounds obtained in Section $1$ of
\cite{ILS}. We consider weight $k$ cuspidal newforms of prime
level $N$ and odd functional equation. We use Theorem
\ref{thm:extended moments} with $n=2$ and
\be\label{eq:phihphibounds} \phi(x) \ = \ \left(\frac{\sin \pi
\sigma x}{\pi \sigma x}\right)^2, \ \ \ \twocase{\hphi(y)\ = \
}{\frac1{\sigma}-\frac{|y|}{\sigma}}{if $|y| <
\sigma$}{0}{otherwise.} \ee Arguing as in Section 1 of \cite{ILS},
we find that as $N\to\infty$ the probability that the zero at the
central point is of order exactly one is at least
$\frac{23}{27}-\gep \approx .8519-\gep$, which is worse than the
lower bound of $\frac{15}{16}-\gep = .9375-\gep$ of \cite{ILS}. To
see this, take $\sigma = 1$ in \eqref{eq:phihphibounds}; while we
need $\sigma < 1$, we can take the limit as $\sigma$ approaches
$1$ (or better yet, introduce a factor of $\gep$ in the
arguments). Let $p_r(N)$ be the percent of odd cuspidal newforms
of weight $k$ and prime level $N$ that have exactly $r$ zeros at
the central point. As our forms are odd, only odd values of $r$
are non-zero. We have \be \sum_{j=0}^\infty p_{2j+1}(N) \ = \ 1, \
\ \ \ \sigma_\phi^2 - R_2(\phi) \ \le \ \frac13 + \gep. \ee

Consider the terms $D(f;\phi) - \< D(f;\phi)\>_-$ in Theorem
\ref{thm:extended moments}. As $\phi(0) = \hphi(0) = 1$ and $\phi$
is non-negative, we see that if there are $r \ge 3$ zeros at the
central point, then \be \begin{array}{ccccc}  D(f;\phi) - \<
D(f;\phi)\>_- & \ \ge \ & r \phi(0) - \< D(f;\phi)\>_- & \ \ge \ & 0
\nonumber\\  & \ge & r\phi(0) - \left(\hphi(0) +
\foh\phi(0)+\gep\right) & \ \ge \ & 0 \nonumber\\  & \ = \ & r -
\frac32-\gep & \ \ge \ & 0.\end{array} \ee If there is exactly one
zero at the central point then $\< D(f;\phi)\>_-$ might exceed
$D(f;\phi)$, and the difference could be negative; if the difference
is negative, when we square we could reverse the inequality. Thus
\bea \frac13 + \gep & \ \ge \ & \sigma_\phi^2 - R_2(\phi)
\nonumber\\ & \ge & \sum_{j=0}^\infty
p_{2j+1}(N) \left(D(f;\phi) - \< D(f;\phi)\>_-\right)^2 \nonumber\\
& \ge & \sum_{j=1}^\infty p_{2j+1}(N) \left(2j+1 -
\frac32-\gep\right)^2 \ \ge \ \left(\frac94-\gep'\right)
\sum_{j=1}^\infty p_{2j+1}(N). \eea Therefore $\sum_{j \ge 1}
p_{2j+1}(N) \le \frac4{27}+\gep''$, or $p_1(N) \ge
\frac{23}{27}-\gep''$.

Our results are better as the order of the zeros increase. A
similar analysis shows the probability that the zero at the
central point is of order at least $5$ is at most
$\frac{4}{147}+\gep \approx .02721+\gep$, which is better than the
upper bound of $\frac1{32}+\gep = .03125+\gep$ implicit in
\cite{ILS}.

\begin{rek}\label{rek:HoffLocknoharm}
In order to obtain bounds on the order of vanishing at the central
point, it is necessary to weigh each cusp form equally (by $c_k
N^{-1}$). While the harmonic weights $\omega_N(f) =
\frac{\Gamma(k-1)}{(4\pi)^{k-1} (f,f)_N}$ are almost constant, by
\cite{I1,HL} they can fluctuate within the family as \be
N^{-1-\gep}\ \ll_k \ \omega_N(f) \ \ll_k \ N^{-1+\gep}; \ee if we
allow ineffective constants we can replace $N^\gep$ with $\log N$
for $N$ large. The difficulty with using harmonic weights is that
the larger weights could all be associated to $f$'s with large (or
small) vanishing at the central point. This is one of the main
reasons we chose to use uniform weights; see also Remark
\ref{rek:whynoharm}.
\end{rek}


\appendix


\section{Handling the Complementary Sum}
\begin{lem}\label{lemcompsumsmall}
Assume GRH for $L(s,f)$ for $f \in H_k^\ast(1)\cup H_k^\ast(N)$. In
the notation of Lemma \ref{lemcomplsum}, if $W=1$ or $N$, $X = N-1$
or $N^\gep$, and $Y = N^\gep$, then
\begin{align}\label{eq:compsumlemmaapp}
& \frac{1}{|H_k^\pm(N)|} \sum_{\substack{q_1 \notdiv N , \dots, q_\ell \notdiv N \\
q_j {\rm distinct}}} \prod_{j=1}^\ell \hphi\left(\frac{\log
q_j}{\log R}\right)^{n_j} \left(\frac{2 \log q_j}{\sqrt{q_j} \log R}
\right)^{n_j}\ \Delta_{k,N}^\infty(W q_1^{m_1} \cdots
q_\ell^{m_\ell}) \nonumber\\ & \ \ \ \ll \ O(N^{-\gep''}W^{-1})
\end{align}
for some $\gep'' > 0$.
\end{lem}

\begin{proof}
From \eqref{eq:number of terms in hkpm} we have that $H_k^\pm(N)
\sim \frac{(k-1)N}{24}$. It suffices to show
\begin{equation}\label{eq:A2sumqcond}
\sum_{(q,N) = 1} \gl_f(q) a_q \ \ll \ (kN)^{\gep'},
\end{equation}
where
\begin{equation}
a_q \ = \
\begin{cases}
\prod_{j=1}^\ell\left(\hphi\left(\frac{\log q_j}{\log R}\right)
\frac{\log q_j}{\sqrt{q_j} \log R}\right)^{n_j} & q =
q_1^{m_1}\cdots q_\ell^{m_\ell},
\text{$q_j \leq R^\ga$ distinct primes, $q_j \notdiv N$}\\
0 & \text{otherwise.}
\end{cases}
\end{equation}
This is because if \eqref{eq:A2sumqcond} holds, Lemma
\ref{lemcomplsum} implies that \be \frac{1}{|H_k^\pm(N)|}
\sum_{(q,WN) = 1 \atop \log q \ll \log N} \Delta_{k,N}^\infty(Wq)
a_q \ \ll \ N^{-1} \cdot N^{1-\gep''} W^{-1} \ \ll \ N^{-\gep''}
W^{-1}. \ee

Without loss of generality, we may relabel so that $q_1 > \dots
> q_\ell$. Up to combinatorial factors, our
sum \eqref{eq:compsumlemmaapp} becomes
\begin{equation}
\sum_{q_\ell = 2}^{R^\alpha} \gl_f(q_\ell^{m_\ell})
\left(\hphi\left(\frac{\log q_\ell}{\log R}\right)\frac{\log
q_\ell}{\sqrt{q_\ell}\log R}\right)^{n_\ell} \sum_{q_{\ell-1} =
q_{\ell}+1}^{R^\ga} \cdots  \sum_{q_1 = q_2+1}^{R^\ga}
\gl_f(q_1^{m_1}) \left(\hphi\left(\frac{\log q_1}{\log
R}\right)\frac{\log q_1}{\sqrt{q_1} \log R}\right)^{n_1}.
\end{equation}
The Generalized Riemann Hypothesis for $L(s,f)$ yields
\begin{equation}
\sum_{p \le P} \frac{\gl_f(p) \log p}{\sqrt{p}} \ \ll \
(kN)^{\gep'/n}
\end{equation}
if $\log P \ll \log kN$ (see equations 2.65--2.66 of \cite{ILS}).
Therefore by partial summation \begin{equation} \sum_{p \le P}
 \hphi\left(\frac{\log p}{\log
R}\right)\frac{\gl_f(p) \log p}{\sqrt{p}} \ \ll \ (kN)^{\gep'/n}
\end{equation} Thus all the sums with $n_j = 1$ are
$\ll (kN)^{\gep'}$. For factors with $n_j
> 1$, the Ramanujan bound for $\gl_f(p)$ gives $|\gl_f(q_j^{m_j})|
\le \tau(q_j^{m_j}) \le n+1$ (as $m_j \le n$), and these prime
sums are then at most
\begin{equation}
\sum_{p\leq P} \frac{\log^{n_j} p}{p \log^{n_j} R}\ \ll \
\frac{\log^{n_j+1} P}{\log^{n_j} R}\ \ll \ (kN)^{\gep'/n}.
\end{equation}
Thus $\sum_q \gl_f(q)a_q \ll (kN)^{\gep'}$, and by Lemma
\ref{lemcomplsum} (and the remarks immediately following it), this
completes the proof.
\end{proof}


\section{Handling the Error Terms in the Moment
Expansion}\label{sec:handerrmomexp}

In order to prove \eqref{eq:centered mmts of D in terms of P}, we
must show that if $\supp(\hphi) \in [-1,1]$ then
\begin{equation}\label{eq:apperrmom}
\left\langle\left(-P(f;\phi) + \ O\left( \frac{\log \log R}{\log
R}\right)\right)^n\right\rangle_\sigma\ =\
(-1)^n\<P(f;\phi)^n\>_\sigma \ + \ O\left( \frac{\log \log R}{\log
R} \right),
\end{equation}
where $\sigma \in \{+,-,\ast\}$ and where $n$ is an integer $\geq
1$. Note the $O$-term on the left hand side of
\eqref{eq:apperrmom} is independent of $f$. Let $P = P(f;\phi)$
and $E = O\left( \frac{\log \log R}{\log R} \right)$. Assume we
know that if $\supp(\hphi) \in [-1,1]$ then
\begin{equation}\label{eq:appmomp2m}
\<P^{2m}\>_\sigma\ =\ O(1)
\end{equation}
for all $0\leq m \le \frac{n}2$. In general, in investigations of
the $n$\textsuperscript{th} centered moments one has $m <
\frac{n}2$ by induction, and handling $m = \frac{n}2$ is possible
-- in fact, this is the expected main term that we evaluate in
\S\ref{sec:NTncentmom}. Expanding, we find
\begin{equation}
\left\langle\left(-P + E\right)^n\right\rangle_\sigma \ = \
(-1)^n\<P^n\>_\sigma \ + \ \sum_{j=1}^n \binom{n}{j} \left\langle
(-P)^{n-j} E^j \right\rangle_\sigma,
\end{equation}
where $E=O\left( \frac{\log \log R}{\log R}\right)$ is independent
of $f$. We show for all $j=1,\dots, n$ that
\begin{equation}
\left\langle (-P)^{n-j} E^j \right\rangle_\sigma\ =\ O(E^j).
\end{equation}
If $n-j$ is even then
\begin{equation}
\left\langle (-P)^{n-j} E^j \right\rangle_\sigma\ \leq\  \<
P^{n-j} \>_\sigma O(E^j)\ =\ O(E^j)
\end{equation}
since we assumed that $\< P^{n-j} \>_\sigma = O(1)$ and that $E$
is independent of $f$ (so it can be taken out of the average). If
$n-j$ is odd, we use the following form of  H\"older's inequality:
if $f,g,\mu$ are positive functions then for $0<\theta<1$,
\begin{equation}
\int f(x) g(x) \mu(x)\;\d x\ \leq\ \left(\int f(x)^{1/\theta}
\mu(x)\;\d x \right)^{\theta} \left(\int g(x)^{1/(1-\theta)}
\mu(x)\;\d x \right)^{1-\theta}.
\end{equation}
Hence
\begin{align}
\left\langle (-P)^{n-j} E^j \right\rangle_\sigma &\ \leq\
\left\langle
|P|^{n-j} E^j \right\rangle_\sigma \nonumber\\
&\ \leq\ \< |P|^{(n-j)/\theta} \>_\sigma^\theta \<
E^{j/(1-\theta)}
\>_\sigma^{1-\theta}\nonumber\\
&\ =\ \< |P|^{(n-j)/\theta} \>_\sigma^\theta E^j.
\end{align}
Now choose $\theta = (n-j)/(n-j+1)<1$, which means $(n-j) / \theta =
n-j+1$. This will be even since $n-j$ is odd, and is clearly less
than or equal to $n$ (as $j \ge 1$). Hence
\begin{equation}
\left\langle (-P)^{n-j} E^j \right\rangle_\sigma\ \leq\ \<
P^{n-j+1} \>_\sigma^{(n-j)/(n-j+1)} E^j\ =\ O(E^j)
\end{equation}
since we assumed that $\< |P|^{n-j+1} \>_\sigma = O(1)$. This
completes the proof of \eqref{eq:centered mmts of D in terms of
P}.


\section{Kloosterman Sum Expansion}\label{sec:KloostSumExpAppC}

As remarked in \cite{ILS}, it is advantageous to employ characters
to a smaller modulus (to modulus $b$ rather than $Nb$) in expanding
the Kloosterman terms.

\begin{lem}\label{lem:usingsmallermodsimp} If $(P,b)=1$ and $(N,b)=1$, then
\be S(m^2,P N; Nb) \ = \ \frac{-1}{\varphi(b)} \sum_{\chi \smod{b}}
\chi(N) G_{\chi}(m^2) G_{\chi}(1) \overline{\chi}(P). \ee
\end{lem}

\begin{proof}
By the orthogonality relation for characters, since $(P,b)=1$ and
$S(m^2,PN;Nb)$ is periodic in $P$ modulo $b$, we may write
\begin{align}
S(m^2,P N; Nb) \ &= \ \frac{1}{\varphi(b)} \sum_{\chi \smod{b}}\
\sideset{}{^*}\sum_{a \bmod b} \chi(a) S(m^2,aN;Nb)
\overline{\chi}(P) \nonumber\\
&=\ \frac{1}{\varphi(b)} \sum_{\chi \smod{b}}\
\sideset{}{^*}\sum_{a \bmod b} \chi(a) \sideset{}{^*}\sum_{d\bmod
Nb}
e(m^2d/Nb) e(aN\overline{d}/Nb) \overline{\chi}(P) \nonumber\\
&= \ \frac{1}{\varphi(b)} \sum_{\chi \smod{b}}\
\sideset{}{^*}\sum_{d \bmod Nb} \chi(d) e(m^2d/Nb)
\sideset{}{^*}\sum_{a \bmod b} \chi(a) e(a/b) \overline{\chi}(P)
\nonumber\\
&= \ \frac{1}{\varphi(b)} \sum_{\chi
\smod{b}}\ \sideset{}{^*}\sum_{d \bmod Nb} \chi(d) e(m^2d/Nb)
G_{\chi}(1) \overline{\chi}(P).
\end{align}
Since $(N,b)=1$ we may replace the sum over $d$ relatively prime to
$Nb$ with $d = u_1N + u_2b$, with $u_1 \bmod b$ relatively prime to
$b$ and $u_2 \bmod N$ relatively prime to $N$. As $\chi$ is a
character modulo $b$, $\chi(u_1N+u_2b) = \chi(u_1N)$. Thus \bea
\sideset{}{^*}\sum_{d\bmod Nb} \chi(d) e(m^2d/Nb) & \ = \ &
\sideset{}{^*}\sum_{u_1 \bmod b} \chi(u_1N) e(m^2u_1/b)
\sideset{}{^*}\sum_{u_2 \bmod N} e(m^2u_2/N) \nonumber\\ & = &
\chi(N) G_{\chi}(m^2) \cdot \left[-1 + \sum_{u_2=0}^{N-1}
e(m^2u_2/N)\right] \nonumber\\ & = & -\chi(N) G_{\chi}(m^2), \eea
because the $u_1$-sum is $G_{\chi}(m^2)$ and $N$ is prime.
Substituting back yields the lemma.
\end{proof}

The reason for using characters with smaller moduli is that we
obtain a savings in estimating the contributions from the
non-principal characters.
\begin{lem}\label{lem:csboundchar} We have
\begin{equation}\label{eq:csboundchar}
\frac1{\varphi(b)}\sum_{\chi \smod{b}} \left|G_\chi(m^2)
G_\chi(1)\right| \ \ll \ \varphi(b) \ \ll \ b;
\end{equation}
\end{lem}

\begin{proof}
{From} the orthogonality of the characters we have
\begin{equation}
\sum_{\chi \smod{b}}\left|G_\chi(n)\right|^2\ =\ \varphi(b)^2,
\end{equation}
and \eqref{eq:csboundchar}  follows from this bound and the
Cauchy-Schwartz inequality.
\end{proof}

Note that if we used characters of modulus $Nb$ then the bound $b$
would be replaced with $Nb$.

\begin{lem}\label{lem:usingsmallermodsimpr}
Assume $(Q,N)= (N,b)=1$, and set $r = (Q,b)$, $b' = b/r$ and $Q' =
Q/r$. If additionally $(r, b') = 1$ then \be S(m^2,Q; Nb) \ = \
\frac{1}{\varphi(Nb/r)} \sum_{\chi \smod{Nb/r}}\overline{\chi}(Q/r)
\chi(r) R(m^2,r) G_{\chi}(m^2) G_{\chi}(1). \ee
\end{lem}

\begin{proof} From the definition of $r$ we have $(Q', b') = 1$ (if
not, $r$ is not the greatest common divisor of $Q$ and $b$). By the
orthogonality relation for characters, since $(Q',Nb')=1$ we may
write
\begin{align}
S(m^2, Q; Nb) \ &= \ S(m^2, Q'r; Nb'r) \nonumber\\ &=\
\frac{1}{\varphi(Nb')} \sum_{\chi \smod{Nb'}}\
\sideset{}{^*}\sum_{a \bmod Nb'} \chi(a)\overline{\chi}(Q')
S(m^2,ar;Nb'r)
\nonumber\\
&=\ \frac{1}{\varphi(Nb')} \sum_{\chi \smod{Nb'}}\
\sideset{}{^*}\sum_{a \bmod Nb'} \chi(a) \overline{\chi}(Q')
\sideset{}{^*}\sum_{d\bmod Nb'r}
e(m^2d/Nb'r) e(ar\overline{d}/Nb'r)  \nonumber\\
&= \ \frac{1}{\varphi(Nb')} \sum_{\chi \smod{Nb'}}\
\sideset{}{^*}\sum_{d \bmod Nb'r} \overline{\chi}(Q') \chi(d)
e(m^2d/Nb'r) \sideset{}{^*}\sum_{a \bmod Nb'} \chi(a) e(a/Nb')
\nonumber\\
&= \ \frac{1}{\varphi(Nb')} \sum_{\chi \smod{Nb'}}\
\sideset{}{^*}\sum_{d \bmod Nb'r}\overline{\chi}(Q') \chi(d)
e(m^2d/Nb'r) G_{\chi}(1).
\end{align}
As $r|Q$ and $(Q,N) = 1$, $(r,N)=1$. Thus $(Nb',r) = 1$, and we may
replace the sum over $d$ relatively prime to $Nb'r$ with $d = u_1Nb'
+ u_2r$, with $u_1 \bmod r$ relatively prime to $r$ and $u_2 \bmod
Nb'$ relatively prime to $Nb'$. As $\chi$ is a character modulo
$Nb'$, $\chi(u_1Nb'+u_2r) = \chi(u_2r)$. Thus \bea
\sideset{}{^*}\sum_{d\bmod Nb'r} \chi(d) e(m^2d/Nb'r) & \ = \ &
\sideset{}{^*}\sum_{u_1 \bmod r} e(m^2u_1Nb'/Nb'r)
\sideset{}{^*}\sum_{u_2 \bmod Nb'} \chi(u_2r) e(m^2u_2r/Nb'r) \nonumber\\
& = & \sideset{}{^*}\sum_{u_1 \bmod r} e(m^2u_1/r) \cdot \chi(r)
\cdot \sideset{}{^*}\sum_{u_2 \bmod Nb'} \chi(u_2)
e(m^2u_2/Nb')\nonumber\\ & = & \chi(r) R(m^2,r) G_{\chi}(m^2), \eea
because by \eqref{eq:defn Rnew} the $u_1$-sum is $R(m^2,r)$ and by
\eqref{eq:gausssum} the $u_2$-sum is $G_{\chi}(m^2)$. Substituting
back yields the lemma.
\end{proof}


\section{Fourier Transform Identities}

Let $\I_{\{|u|\le 1\}}$ be the characteristic function of
$[-1,1]$. Let $S(x) = \frac{\sin \pi x}{\pi x}$. Note that
\begin{equation}\label{eq:S(2x) as Fourier transform}
S(2x)\ =\ \intinf \foh \I_{\{|u|\le 1\}}\ e^{2\pi\i x u} \;\d u,
\end{equation}
so $S(2x)$ and $\foh \I_{\{|u|\leq 1\}}$ are a Fourier transform
pair. All test functions below will be even Schwartz functions
whose Fourier transforms have finite support. We have made much
use of a certain Fourier transform identity; we give the proof
here for completeness.
\begin{lem}\label{thm:neededn=2int} We have
\begin{multline}\label{eq:FT id1}
\int_{-\infty}^\infty \phi(x)^n S(2x)
\;\d x - \frac12 \phi(0)^n \\
=\ -\frac12\idotsint_{-\infty}^\infty \hphi(u_1) \cdots \hphi(u_n)
\I_{\{|u_1+\cdots+u_n| > 1\}} \;\d u_1 \cdots \;\d u_n .
\end{multline}
\end{lem}

\begin{proof}
This lemma follows from Plancherel's identity, which states that
if $f$ and $g$ are Schwartz functions (in fact it is true for a
much larger class of functions) then
\begin{equation}
\int f(x) g(x) \;\d x\ =\ \int \widehat f(u) \widehat g(u) \;\d u .
\end{equation}
In this particular case it is more complicated since we are
integrating $n+1$ functions. We obtain
\begin{equation}\label{eq:appC eq 1}
\phi(0)^n\ =\ \idotsint_{-\infty}^\infty \hphi(u_1) \cdots
\hphi(u_n) \;\d u_1 \cdots \;\d u_n
\end{equation}
and
\begin{equation}\label{eq:appC eq 2}
\int_{-\infty}^\infty \phi(x)^n \frac{\sin 2\pi x}{2 \pi x} \;\d x \
=\ \frac12\idotsint_{-\infty}^\infty \hphi(u_1) \cdots \hphi(u_n)
\I_{\{|u_1+\cdots+u_n|\leq 1\}} \;\d u_1 \cdots \;\d u_n,
\end{equation}
where we have used \eqref{eq:S(2x) as Fourier transform} and
(repeatedly)
\begin{equation}
\widehat{f g}(u)\ =\ \intinf \widehat f(v) \widehat g(u-v) \;\d v
.
\end{equation}
Combining \eqref{eq:appC eq 1} and \eqref{eq:appC eq 2} yields
\begin{multline}
\int_{-\infty}^\infty \phi(x)^n S(2x) \;\d x - \frac12 \phi(0)^n
\\
=\ \frac12\idotsint_{-\infty}^\infty \hphi(u_1) \cdots \hphi(u_n)
\left( \I_{\{|u_1+\cdots+u_n|\leq 1\}}-1\right) \;\d u_1 \cdots \;\d
u_n\\
=\ -\frac12\idotsint_{-\infty}^\infty \hphi(u_1) \cdots \hphi(u_n)
\I_{\{|u_1+\cdots+u_n|> 1\}} \;\d u_1 \cdots \;\d u_n,
\end{multline}
which is \eqref{eq:FT id1}.
\end{proof}


\section{Increasing the Support in Theorem \ref{thm:nosplitmoments}}

As it stands, Theorem \ref{thm:nosplitmoments} holds for
$\supp(\hphi) \subset \left(-\frac2n\left(1 - \frac1{2k}\right),
\frac2n\left(1 - \frac1{2k}\right)\right)$. We show how a more
careful book-keeping and using GRH for Dirichlet $L$-functions
allows us to remove the factors of $\frac1{2k}$ for $n \ge 2$ and
$2k \ge n$. In particular, we prove

\begin{thm}\label{thm:nosplitbetterk} Assume GRH for Dirichlet
$L$-functions. If $2k \ge n$ then Theorem \ref{thm:nosplitmoments}
holds for even Schwartz test functions supported in $(-\frac2n,
\frac2n)$.
\end{thm}

In proving Lemma \ref{lem:S1} (which is equivalent to Theorem
\ref{thm:nosplitmoments}) we showed, without any restriction on the
support of $\hphi$, that for $S_1^{(n)}$ defined as in
\eqref{eq:S_1}, then under GRH for $L(s,f)$
\begin{equation}
\lim_{\substack{N\to\infty \\ N \text{ prime}}} S_1^{(n)} =
\begin{cases}
\frac{(2m)!}{2^m m!}\left( \ 2\int_{-\infty}^{\infty} |y| \hphi(y)^2
\;\d
y \right)^{2m} + E(n) & \text{if $n=2m$ is even}\\
E(n) & \text{if $n$ is odd,}
\end{cases}
\end{equation}
where $E(n)$ is made up of a linear combination of terms like
\begin{multline}\label{eq:enexpmn}
\lim_{\substack{N\to\infty \\ N \text{ prime}}}
\sum_{\substack{q_1, \dots, q_\ell \\
q_j  \text{ distinct primes}}} \prod_{j=1}^\ell
\hphi\left(\frac{\log q_j}{\log R}\right)^{n_j} \left(\frac{2 \log
q_j}{\sqrt{q_j} \log R}
\right)^{n_j} \\
\times \sum_{m \le N^\gep} \frac{2\pi \i^k}{m} \sum_{b=1}^\infty
\frac{S(m^2,q_1^{m_1}\cdots q_\ell^{m_\ell},Nb)}{Nb}
J_{k-1}\left(4\pi m\frac{\sqrt{q_1^{m_1}\cdots
q_\ell^{m_\ell}}}{Nb}\right),
\end{multline}
where $n_j\geq 1$ with $n_1+\cdots+n_\ell = n$, and $m_j \leq n_j$
with $m_j \equiv n_j \bmod 2$. Lemma~\ref{lem:S1} followed from
this by showing $E(n) = 0$ if $\supp(\hphi) \subset
(-\frac{2k-1}{kn},\frac{2k-1}{kn})$, via the bound
\eqref{eq:estimate Kloosterman} on the Kloosterman sum, and the
bound from Lemma~\ref{lem:Bessel}(\ref{lb:3}) on the Bessel
function. We prove Theorem \ref{thm:nosplitbetterk} by showing, in
a sequence of lemmas, that whenever $2k \ge n$ then $E(n) = 0$ for
$\supp(\hphi) \subset (-\frac2n, \frac2n)$.

For simplicity we write $Q = q_1^{m_1} \dots q_\ell^{m_\ell}$ and
$\supp(\hphi) \subset [-\sigma, \sigma] \subset (-\frac2n,
\frac2n)$. As $n \ge 2$, $(Q,N) = 1$. Set $r = (Q,Nb)$, $Q' = Q/r$
and $b' = b/r$. If additionally $(r, b/r) = 1$ then Lemma
\ref{lem:usingsmallermodsimpr} yields \be\label{eq:appKloosttouse}
S\left(m^2,\frac{Q}{r} r, \frac{Nb}{r} r\right) \ = \
\frac{1}{\varphi(Nb/r)} \sum_{\chi \smod{Nb/r}} \overline{\chi}(Q/r)
\chi(r) R(m^2,r) G_{\chi}(m^2) G_{\chi}(1). \ee

We sketch the proof of Theorem \ref{thm:nosplitbetterk}. In
\S\ref{lem:enrbrterm} we handle the terms in $E(n)$ with $(r, b/r) >
1$ (and thus the expansion of Lemma \ref{lem:usingsmallermodsimpr}
is unavailable), thereby reducing the proof to an analysis of the
terms with $(r, b/r) = 1$. In \S\ref{sec:appErestrictb} we show we
may truncate the $b$-sum at $N$; this is useful as some later terms
are $\sum_b b^{-1}$. We then show in \S\ref{sec:appEkloost} that we
may assume $r=1$, and then use Lemma \ref{lem:usingsmallermodsimpr}
to expand the Kloosterman sums. The proof is completed in
\S\ref{sec:appEGRH} where we bound the contributions from the
character sums arising from the Kloosterman expansions; it is here
that we must assume GRH for Dirichlet $L$-functions.

\subsection{Bounding the terms with $(r, b/r) > 1$}

\begin{lem}\label{lem:enrbrterm}
Notation as above, the contribution to $E(n)$ from terms with $(r,
b/r) > 1$ is negligible for $\supp(\hphi) \subset (-\frac2n,
\frac2n)$, provided that $8k - 2 \ge n$. \end{lem}

\begin{proof} As $Q = q_1^{m_1} \cdots q_\ell^{m_\ell}$ is a product
over distinct primes, if $(r, b/r) > 1$ then the square of some
$q_j$ divides $b$. Without loss of generality we may assume $b =
q_1^2 v$. For $\supp(\hphi)$ and $k$ as in the lemma, we show the
contribution to $E(n)$ (equation \eqref{eq:enexpmn}) from each tuple
$(n_1, \dots, m_\ell)$ is negligible. This proves the lemma as the
number of such tuples depends only on $n$ and not $N$.

We use $J_{k-1}(x) \ll x^{k-1}$. This is the best available Bessel
bound for our purposes, as our argument is at most $N^\gep$. We use
\eqref{eq:estimate Kloosterman} to bound the Kloosterman sum. Thus
the contribution from such a tuple is \bea
E(\overrightarrow{n},\overrightarrow{m}) &\ll\ & \sum_{m \le N^\gep}
\frac1m \sum_{v=1}^\infty \frac{(Nq_1^2v / q_1)^{\foh+\gep}}{Nq_1^2
v} \sum_{q_1, \dots, q_\ell = 2}^{N^\sigma} (q_1 \cdots
q_\ell)^{-\frac{n_j}2} \left( \frac{m\sqrt{q_1 \cdots
q_\ell}}{Nq_1^2v}\right)^{k-1} \nonumber\\ & \ll & N^{\foh-k+\gep'}
\sum_{m \le N^\gep} m^{k-2} \sum_{v=1}^\infty
\frac{1}{v^{k-\frac12}} \sum_{q_1=2}^{N^\sigma}
q_1^{-\frac{n_1}2+\frac{m_1}{2}(k-1) -2k + \foh} \prod_{j=2}^\ell
\sum_{q_j=2}^{N^\sigma}
q_j^{-\frac{n_j}2+\frac{m_j}{2}(k-1)}.\nonumber\\ \eea The worst
case is when $m_j = n_j$. As the $v$-sum is $O(1)$, the $m$-sum is
$O(N^{(k-1)\gep})$, and $\sum n_j = n$, \bea
E(\overrightarrow{n},\overrightarrow{m}) & \ll\ & N^{\foh-k+\gep''}
N^{(\frac{nk}2 - n + \ell - 2k + \foh)\sigma}. \eea The worst case
is when $\ell = n$, and we find \be
E(\overrightarrow{n},\overrightarrow{m}) \ \ll \ N^{\foh\left((nk -
4k + 1)\sigma - (2k-1)\right) + \gep'''}. \ee This is negligible
provided that $\sigma < \frac{2k-1}{nk-4k+1}$. If $8k - 2 \ge n$,
simple algebra shows $\frac{2k-1}{nk-4k+1} \ge \frac2n$.
\end{proof}

\subsection{Restricting the $b$-sum}\label{sec:appErestrictb}

\begin{lem}\label{lem:appEasize}
If $a \ge \frac1{2k-3}$, then the contribution to $E(n)$ (equation
\eqref{eq:enexpmn}) from the $b \ge N^a$ terms is negligible for any
admissible $(m,n)$-tuple. \end{lem}

\begin{proof} Using $J_{k-1}(x) \ll x^{k-1}$ and $S(m^2,Q,bN) \ll (bN)^\foh N^\gep$,
the contribution to \eqref{eq:enexpmn} from terms with $b \ge N^a$
is \bea & \ll \ & N^{-1} \sum_{m \le N^\gep} m^{k-2}
\prod_{j=1}^\ell \left( \sum_{q_j =2}^{N^\sigma}
q_j^{\frac{m_j}2(k-1) - \frac{n_j}2}\right) \sum_{b=N^a}^\infty
\frac{b^\foh N^{\foh + \gep}}{b \cdot b^{k-1}} \frac1{N^{k-1}}
\nonumber\\ & \ll \ & N^{-k+\foh+\gep'} N^{\frac{nk\sigma}2}
N^{\left(-k+\frac32\right)a}, \eea because the worst case is when
$\ell=n$ and $m_j=n_j$. If $\sigma < \frac2n$ then for $\gep'$
sufficiently small there is no contribution provided \be \foh +
\left(-k+\frac32\right)a \ \le \ 0, \ee which means $a \ge
\frac{1}{2k-3}$. As $k\ge 2$, $a \le 1$.
\end{proof}

\subsection{Restricting the $r$-sum and Expanding the Kloosterman Sum}\label{sec:appEkloost}

The proof of Theorem \ref{thm:nosplitbetterk} is therefore reduced
to showing that there is no contribution from admissible
$(m,n)$-tuples as $N\to\infty$ in \bea\label{eq:s1toboundfixedm2}
& & E(\overrightarrow{n},\overrightarrow{m}) \ = \ \nonumber\\
& & \sum_{\substack{q_1, \dots, q_\ell \\
q_j {\rm distinct}}} \prod_{j=1}^\ell \hphi\left(\frac{\log
q_j}{\log R}\right)^{n_j} \left(\frac{2 \log q_j}{\sqrt{q_j} \log R}
\right)^{n_j} \sum_{m \le N^\gep} \frac{2\pi \i^k}{m} \sum_{b= 1
\atop (Q,b) = r, (r,b/r)=1}^{N^{1/(2k-3)}} \frac{S(m^2,Q,Nb)}{bN}
J_{k-1}\left(4\pi m \frac{\sqrt{Q}}{bN}\right);\nonumber\\ \eea by
an admissible $(m,n)$-tuple we mean $\sum n_j = n$, $m_j \le n_j$
and $m_j \equiv n_j \bmod 2$. As $r = (Q,b)$ and $Q = q_1^{m_1}
\cdots q_\ell^{m_\ell}$, we may write \be r \ = \ q_1^{c_1} \cdots
q_\ell^{c_\ell}, \ \ \ c_j \in \{0,\dots, m_j\}. \ee For a given
$m$-tuple $(m_1,\dots,m_\ell)$, the number of $c$-tuples is
$\prod_{j=1}^\ell (m_j+1) \ll (2n)^n$. We show

\begin{lem} Notation as above, a $c$-tuple with $\sum c_j = c$ has a
negligible contribution to
$E(\overrightarrow{n},\overrightarrow{m})$ for $\supp(\hphi)
\subset(-\frac2n, \frac2n)$, provided that $4ck \ge n$. In
particular, if $4k \ge n$ then to bound
$E(\overrightarrow{n},\overrightarrow{m})$ it suffices to consider
only the contributions from the $c$-tuple $(0,\dots,0)$ (i.e.,
$r=1$).
\end{lem}

\begin{proof} We use \eqref{eq:estimate Kloosterman} to bound the
Kloosterman sum, and $J_{k-1}(x) \ll x^{k-1}$. For a fixed $c$-tuple
with $\sum c_j = c$ we have $b = b'r$. We insert absolute values and
ignore the condition $(r,b/r) = 1$, as this only increases the sums
below. We have \bea
E(\overrightarrow{n},\overrightarrow{m},\overrightarrow{c}) & \ \ll
\ & \sum_{m \le N^\gep} \frac1m \sum_{q_1,\dots,q_\ell =
2}^{N^\sigma} \sum_{b'=1}^{N} \frac{(Nb')^{\foh+\gep}}{Nb'
q_1^{c_1}\cdots q_\ell^{c_\ell}} \prod_{j=1}^\ell q_j^{-\frac{n_j}2}
\cdot \left(\frac{m \sqrt{q_1^{m_1} \cdots
q_\ell^{m_\ell}}}{b'q_1^{c_1}\cdots
q_\ell^{c_\ell}N}\right)^{k-1}.\nonumber\\ \eea The $m$-sum is
$O(N^{(k-1)\gep})$, the $b'$-sum is $O(1)$, and the worst case is
when each $m_j = n_j$. Thus \bea
E(\overrightarrow{n},\overrightarrow{m},\overrightarrow{c}) & \ \ll
\ & N^{\foh-k+\gep'} \prod_{j=1}^\ell \sum_{q_j=2}^{N^\sigma}
q_j^{\frac{n_jk}2-n_j - c_jk} \nonumber\\ & \ll & N^{\foh-k+\gep'}
N^{(\frac{nk}2-n+\ell-ck)\sigma}. \eea As usual, the worst case is
when $\ell = n$, and we find \be
E(\overrightarrow{n},\overrightarrow{m},\overrightarrow{c}) \ \ll \
N^{\foh\left((nk-2ck)\sigma - (2k-1)\right) + \gep''}; \ee this is
negligible provided that $\sigma < \frac{2k-1}{nk-2ck}$. If $4ck \ge
n$ then $\frac{2k-1}{nk-2ck} \ge \frac2n$. Thus all $c$-tuples with
$\sum c_j \ge 1$ yield negligible contributions for $\supp(\hphi)
\subset (-\frac2n, \frac2n)$, provided that $4k \ge n$. \end{proof}

\begin{rek} While assuming $4k \ge n$ is more restrictive than $8k \ge
n+2$ (the relation from Lemma \ref{lem:enrbrterm}), this allows us
to take $r=1$ below, and greatly simplifies the arguments. At the
cost of a more involved argument one could analyze the $c$-tuples
where $\sum c_j \in \{1,2\}$. \end{rek}


Thus the proof of Theorem \ref{thm:nosplitbetterk} is reduced to
bounding
$E(\overrightarrow{n},\overrightarrow{m},\overrightarrow{0})$. In
this case $r=1$, so $Q' = Q$ and $b' = b$. From \eqref{eq:defn Rnew}
we have $R(m^2,1)=1$. By Lemma \ref{lem:usingsmallermodsimpr}, we
are left with bounding \bea\label{eq:s1toboundfixedm3}
E(\overrightarrow{n},\overrightarrow{m},\overrightarrow{0}) & \ = \
&
\sum_{\substack{q_1, \dots, q_\ell \\
q_j {\rm distinct}}} \prod_{j=1}^\ell \hphi\left(\frac{\log
q_j}{\log R}\right)^{n_j} \left(\frac{2 \log q_j}{\sqrt{q_j} \log R}
\right)^{n_j}\sum_{m \le N^\gep} \frac{2\pi \i^k}{m}\nonumber\\ & &
\ \ \ \ \ \cdot\  \sum_{b= 1 \atop (Q,b)=1}^{N^{1/(2k-3)}}
\sum_{\chi \bmod Nb \atop \chi^{m_j} \neq \chi_0}
\frac{\overline{\chi}(Q) G_\chi(m^2)G_\chi(1)}{Nb\varphi(Nb)}
J_{k-1}\left(4\pi m \frac{\sqrt{Q}}{bN}\right).\ \eea

\subsection{Using GRH for Dirichlet
$L$-Functions}\label{sec:appEGRH}

We use GRH for Dirichlet $L$-functions to show
\eqref{eq:s1toboundfixedm3} is negligible in the desired range. Note
\be \overline{\chi}(Q) \ = \ \overline{\chi^{m_1}}(q_1) \cdots
\overline{\chi^{m_\ell}}(q_\ell). \ee There is a slight complication
due to the fact that $\chi$ may not be the principal character, but
a $\chi^{m_j}$ is the principal character. We say $\chi$ is a
\textbf{bad character} if $\chi^{m_j}$ is the principal character
for at least one non-zero $m_j$; otherwise $\chi$ is a \textbf{good
character}. Fortunately, as $(N,b) = 1$ and $N$ is prime, for each
admissible $(m,n)$-tuple the number of bad characters is $O_n(b)$;
this follows from the structure theorem for finite abelian groups
and the fact that our prime $N$ is relatively prime to $b$.

The proof of Theorem \ref{thm:nosplitbetterk} is completed in the
following two lemmas (Lemmas \ref{lem:emnbad} and
\ref{lem:emngood}), which combined show that for an admissible
$(m,n)$-tuple, the contribution to
$E(\overrightarrow{n},\overrightarrow{m},\overrightarrow{0})$ from
both the bad and the good characters is negligible.

\begin{lem}\label{lem:emnbad} For an admissible $(m,n)$-tuple, the contribution to
$E(\overrightarrow{n},\overrightarrow{m},\overrightarrow{0})$ from
the bad characters is negligible for $\supp(\hphi) \subset
(-\frac2n,\frac2n)$.
\end{lem}

\begin{proof} There are at most $O_n(b)$ bad characters. We insert
absolute values in all sums below. Thus the worst case is when each
$\chi^{m_j}$ is the principal character $\chi_0$.  As $G_\chi(a)$ is
a Gauss sum for a character of modulus $bN$, we have $G_\chi(a) \ll
\sqrt{bN}$ if $\chi$ is not the principal character; otherwise by
\eqref{eq:defn Rnew} and $m \le N^\gep$ we find $G_{\chi_0}(m^2)
G_{\chi_0}(1) \ll m^4 \cdot 1 \ll N^{4\gep}$. We use $J_{k-1}(x) \ll
x^{k-1}$ and insert absolute values. Thus we only increase the sum
when we remove the condition that $(Q,b) = 1$. The contribution to
$E(\overrightarrow{n},\overrightarrow{m},\overrightarrow{0})$ from
the bad characters is bounded by \bea\label{eq:s1toboundfixedm4}
E_{\rm
bad}(\overrightarrow{n},\overrightarrow{m},\overrightarrow{0}) & \
\ll\ &\sum_{m \le N^\gep} \frac{m^{k-1}}{m} \sum_{b=
1}^{N^{1/(2k-3)}} \frac1{bN} \sum_{\chi \ {\rm bad}}
\frac{|G_\chi(m^2)G_\chi(1)|}{\varphi(bN)}\nonumber\\ & & \ \ \ \
\cdot\ \prod_{j=1}^\ell \left( \sum_{q_j=2}^{N^\sigma}\chi_0(q_j)
\log q_j \cdot q_j^{\frac{m_j(k-1)}2-\frac{n_j}2} \right) \cdot
\frac1{b^{k-1}N^{k-1}}. \nonumber\\ \eea As usual, the worst case is
when each $m_j=n_j$ and $\ell=n$. The product of the $q_j$-sums is
at most $N^{\frac{nk}2\sigma}$. The sum over the bad characters is
$\ll O(b)$ as each summand
$\frac{|G_\chi(m^2)G_\chi(1)|}{\varphi(bN)}$ is $O(1)$. The $m$-sum
is $O(N^{(k-1)\gep})$, and as $k \ge 2$, the $b$-sum is $\sum_{b \le
N} b^{-1} \ll \log N$. Thus \bea E_{\rm
bad}(\overrightarrow{n},\overrightarrow{m},\overrightarrow{0})\ \ll
\ N^{-k+\gep'} N^{\frac{nk}2\sigma}, \eea which is negligible
provided that $\sigma < \frac2n$.
\end{proof}

\begin{lem}\label{lem:emngood} Assume GRH for all Dirichlet $L$-functions.
For an admissible $(m,n)$-tuple and $2k \ge n$, the good characters'
contribution to
$E(\overrightarrow{n},\overrightarrow{m},\overrightarrow{0})$ is
negligible for $\supp(\hphi) \subset (-\frac2n, \frac2n)$.
\end{lem}

\begin{proof} A modification of the argument in Lemma \ref{lem:csboundchar}
gives \be\label{eq:csboundcharbN}
\frac1{\varphi(bN)}\sum_{\smallsubstack{\chi \smod{bN}}{\chi \ {\rm
good}}} \left|G_\chi(m^2) G_\chi(1)\right| \ \ll \ \varphi(bN) \ \ll
\ bN. \ee  We have \bea\label{eq:s1toboundfixedm6} E_{\rm
good}(\overrightarrow{n},\overrightarrow{m},\overrightarrow{0}) & \
= \ &
\sum_{\substack{q_1, \dots, q_\ell \\
q_j {\rm distinct}}} \prod_{j=1}^\ell \hphi\left(\frac{\log
q_j}{\log R}\right)^{n_j} \left(\frac{2 \log q_j}{\sqrt{q_j} \log R}
\right)^{n_j}\sum_{m \le N^\gep} \frac{2\pi \i^k}{m}\nonumber\\ & &
\ \ \ \ \ \ \ \ \cdot\  \sum_{b= 1 \atop (Q,b)=1}^{N^{1/(2k-3)}}
\sum_{\chi \bmod Nb} \frac{\overline{\chi}(Q)
G_\chi(m^2)G_\chi(1)}{Nb\varphi(Nb)} J_{k-1}\left(4\pi m
\frac{\sqrt{Q}}{bN}\right).\nonumber\\ \eea As each $\chi$ is a
character modulo $b$, if $(Q,b) > 1$ then $\overline{\chi}(Q) = 0$.
Thus we may drop the condition that $(Q,b) = 1$.

We first show that we may assume each $m_j \neq 0$ if $2k \ge n$;
note if an $m_j=0$ then $\chi^{m_j}$ is the principal character for
a trivial reason. Without loss of generality, assume $m_1 = 0$.
Since $m_j \equiv n_j \bmod 2$, $n_1 \ge 2$ and $\ell \le n-1$. We
use $J_{k-1}(x) \ll x^{k-1}$; as $m_1 = 0$ there are no factors of
$q_1$ arising from the argument of the Bessel function. Thus the
$q_1$-sum in $E_{\rm
good}(\overrightarrow{n},\overrightarrow{m},\overrightarrow{0})$ is
$O(\log N)$. The worst case is when each remaining $m_j = n_j$ and
$\ell = n-1$ (thus $n_1 = 2$). By \eqref{eq:csboundcharbN},
$\sum_{\chi \bmod Nb} \frac{\overline{\chi}(Q)
G_\chi(m^2)G_\chi(1)}{Nb\varphi(Nb)} \ll 1$. The $b$-sum is $O(1)$,
the $m$-sum is $O(N^{(k-1)\gep})$, and we have a contribution
bounded by \be \ll \ N^{1-k+\gep'} N^{\frac{(n-2)k}2\sigma} \ = \
N^{\foh\left( (n-2)k\sigma - (2k-2)\right) + \gep''}. \ee This is
negligible for $\sigma < \frac{2k-2}{(n-2)k}$, and if $2k \ge n$
then $\frac{2k-2}{(n-2)k} \ge \frac2n$.

Thus we may now assume each $m_j\neq 0$ and each $\chi^{m_j} \neq
\chi_0$. We fix a $b$ and consider the $q_j$-sums. We use partial
summation and GRH for Dirichlet $L$-functions to convert the
$q_j$-sums to integrals. Since $q_j < N$ and $\chi^{m_j}$ is not the
principal character, under GRH we have for $u_j \le N$ that \be
\sum_{q_j \le u_j} \overline{\chi}(q_j^{m_j}) \ = \ H(u_j) \ \ll \
u_j^{\foh} N^\gep, \ee  where $H(u)$ is a non-differentiable
function. This and the compact support of $\hphi$ imply that \bea &
& \sum_{q_1=2}^{N^\sigma} \overline{\chi}(q_1^{m_1}) \cdot
\hphi\left(\frac{\log q_1}{\log R}\right)^{n_1} \left(\frac{2 \log
q_1}{\sqrt{q_1} \log R} \right)^{n_1} J_{k-1}\left(4\pi m
\frac{\sqrt{q_1^{m_1} q_2^{m_2}\cdots
q_\ell^{m_\ell}}}{bN}\right)\nonumber\\ & &  =\ \int_1^{N^\sigma}
H(u_1) \frac{\d}{\d u_1}\left[ \left( \hphi\left(\frac{\log
u_1}{\log R}\right)\frac{2 \log u_1}{\sqrt{u_1} \log R}
\right)^{n_1} J_{k-1}\left(4\pi m \frac{\sqrt{u_1^{m_1}q_2^{m_2}
\cdots q_\ell^{m_\ell}}}{bN}\right)\right]\d u_1.\nonumber\\ \eea
Proceeding in this manner we find \bea\label{eq:appEtemputrf} & &
\prod_{j=1}^\ell \sum_{q_j=2}^{N^\sigma} \overline{\chi}(q_j^{m_j})
\cdot \hphi\left(\frac{\log q_j}{\log R}\right)^{n_j} \left(\frac{2
\log q_j}{\sqrt{q_j} \log R} \right)^{n_j} J_{k-1}\left(4\pi m
\frac{\sqrt{Q}}{bN}\right) \nonumber\\ & & = \ \int
\cdots\int_1^{N^\sigma} H(u_1)\cdots H(u_\ell) \nonumber\\ & & \ \ \
\cdot \frac{\d^\ell}{\d u_1\cdots \d u_\ell} \prod_{j=1}^\ell \left[
\left(\hphi\left(\frac{\log u_j}{\log R}\right) \frac{2 \log
u_j}{\sqrt{u_j} \log R} \right)^{n_j} J_{k-1}\left(4\pi m
\frac{\sqrt{u_1^{m_1} \cdots u_\ell^{m_\ell}}}{bN}\right)\right]\d
u_1 \cdots \d u_\ell\nonumber\\ & & = \ I(m,b,N;
\overrightarrow{n},\overrightarrow{m}). \eea Therefore
\bea\label{eq:s1toboundfixedm7} E_{\rm
good}(\overrightarrow{n},\overrightarrow{m},\overrightarrow{0}) & \
= \ & \sum_{m \le N^\gep} \frac{|2\pi \i^k|}{m}\sum_{b= 1 \atop
(Q,b)=1}^{N^{1/(2k-3)}} \sum_{\chi \bmod Nb}
\frac{|G_\chi(m^2)G_\chi(1)|}{Nb\varphi(Nb)}\ I(m,b,N;
\overrightarrow{n}, \overrightarrow{m}). \nonumber\\ \eea

As $\sigma < \frac2n$, the derivative in \eqref{eq:appEtemputrf} is
bounded by \be \sum_{t=1}^\ell \left|\prod_{j=1}^\ell
u_j^{-\frac{n_j}2-1} J_{k-1}^{(t)}\left(4\pi m \frac{\sqrt{u_1^{m_1}
\cdots u_\ell^{m_\ell}}}{bN}\right)\right|. \ee This follows because
each derivative of a $\left(\hphi\left(\frac{\log u_j}{\log
R}\right) \frac{2 \log u_j}{\sqrt{u_j} \log R} \right)^{n_j}$ with
respect to $u_j$ decreases the exponent of $u_j$ by $1$. If the
differentiation hits the Bessel piece, we get the derivative of the
Bessel function, as well as a factor of $\frac{4\pi m \cdot
u_1^{m_1/2} \cdots u_\ell^{m_\ell/2}}{bN} \frac1{u_j}$. Because
$\sigma < \frac2n$ and $m \ll N^\gep$, the first factor is at most
$O(N^\gep)$ and thus this factor is bounded by $N^{\gep} u_j^{-1}$.
If additional differentiations with respect to $u_{h}$ hit powers of
$\frac{4\pi m \cdot u_1^{m_1/2} \cdots u_\ell^{m_\ell/2}}{bN}$, the
effect is still just to reduce the exponent of $u_h$ by $1$.
Further, by Lemma \ref{lem:Bessel}(\ref{lb:5}), $2J'_\nu(x) =
J_{\nu-1}(x) - J_{\nu+1}(x)$. Using this relation $\ell$ times, as
well as the bound $J_{k-1}(x) \ll x$, we find that $J_{k-1}^{(t)}(x)
\ll x$ and thus the derivative term in \eqref{eq:appEtemputrf} is
\be \ll \ (u_1\cdots u_\ell)^{-1} \prod_{j=1}^\ell
u_j^{-\frac{n_j}2} \cdot \frac{4\pi m \cdot u_1^{\frac{m_1}2} \cdots
u_\ell^{\frac{m_\ell}2}}{bN} \ \ll \ \frac{m}{bN u_1 \cdots u_\ell},
\ee as the worst case is when $m_j = n_j$. Thus, inserting absolute
values and approximating the Bessel function as above, the integral
in \eqref{eq:appEtemputrf} satisfies \be I(m,b,N;
\overrightarrow{n}, \overrightarrow{m}) \ \ll \ N^{\gep'} \int
\cdots \int_1^{N^\sigma} \frac{m}{bN}\ \frac{\d u_1 \cdots \d
u_\ell}{\sqrt{u_1\cdots u_\ell}} \ \ll \ \frac{m}{bN}\
N^{\frac{\ell\sigma}2+\gep'}.\ee As usual, the worst case is when
$\ell = n$, and $I(m,b,N; \overrightarrow{n}, \overrightarrow{m})
\ll N^{\frac{n\sigma}2-1+\gep'} m/b$. Substituting this bound into
\eqref{eq:s1toboundfixedm7} yields \bea E_{\rm
good}(\overrightarrow{n},\overrightarrow{m}, \overrightarrow{0}) & \
\ll \ & \sum_{m \le N^\gep} \sum_{b= 1 \atop (Q,b)=1}^{N^{1/(2k-3)}}
\frac1b \left[\sum_{\chi \bmod Nb}
\frac{|G_\chi(m^2)G_\chi(1)|}{Nb\varphi(Nb)} \right]
N^{\frac{n\sigma}2-1+\gep'}.  \nonumber\\ \eea

By \eqref{eq:csboundcharbN} the bracketed character sum is $O(1)$.
Thus the $b$-sum is $O(\log N)$ and the $m$-sum is $O(N^\gep)$. Thus
\bea E_{\rm good}(\overrightarrow{n},\overrightarrow{m},
\overrightarrow{0}) & \ \ll \ & N^{\frac{n\sigma}2-1+\gep''}, \eea
which is negligible for $\sigma < \frac2n$.
\end{proof}


\section*{Acknowledgements}
Both authors thank AIM for its generous support; the second named
author also wishes to thank Princeton and Ohio State. This research
was partially written up while the first named author was at the
Isaac Newton Institute, funded by EPSRC Grant N09176.  We would also
like to thank Henryk Iwaniec, Wenzhi Luo, Ze\'ev Rudnick, Peter
Sarnak and Soundararajan for many enlightening discussions, and the
referees for a very thorough reading and numerous helpful comments.


\end{document}